\documentclass[12pt,leqno]{article}
\usepackage{amssymb}
\usepackage{amsmath}

\def\t{{\,}^t}
\def\la{\lambda}
\def\al{\alpha}
\def\om{\omega}
\def\ve{\varepsilon}
\def\vp{\varphi}

\def\dom{{\rm dom}}

\def\sp{{\frak s\frak p}}
\def\spec{{\rm spec}\,}

\def\Sp{{\rm Sp}\,}
\def\Re{{\rm Re}\,}
\def\Im{{\rm Im}\,}
\def\Ker{{\rm Ker\,}}
\def\trans{\,{}^t\!} 
\def\sgn{{\rm sgn}\,}

\def\D{{\cal D}}
\def\E{{\cal E}}
\def\L{{\cal L}}
\def\S{{\cal S}}
\def\L{{\cal L}}

\def\J{{\cal J}}
\def\Q{{\cal Q}}
\def\D{{\cal D}}

\def\h{{\frak h}}

\def\de{\partial}

\def\oS{{\overline{S}}}

\def\CC{{\mathbb  C}}

\def\HH{{\mathbb  H}}

\def\NN{{\mathbb  N}}

\def\RR{{\mathbb  R}}

\def\ZZ{{\mathbb  Z}}
\def\11{{\mathbb  1}}

\def\1{{\bf 1}}

\def \trans{\,{}^t\!}

\def\l{\langle}
\def\r{\rangle}

\def\be{\begin{enumerate}}
\def\ee{\end{enumerate}}
\def\noi{\noindent}

\newtheorem{theorem}{Theorem}[section]
\newtheorem{proposition}[theorem]{Proposition}
\newtheorem{example}[theorem]{Example}

\newtheorem{lemma}[theorem]{Lemma}
\newtheorem{remark}[theorem]{Remark}
\newtheorem{remarks}[theorem]{Remarks}
\newtheorem{corollary}[theorem]{Corollary}

\begin{document}

\title{Solvability of dissipative second order left-invariant differential
operators    on the Heisenberg group}
\author{Detlef M\"uller\\[2mm]}
\date{}
\maketitle

\begin{abstract}
We prove local solvability for large classes of operators of the form
$$
L=\sum_{j,k=1}^{2n}a_{jk}V_jV_k+i\alpha U,
$$
where the $V_j$ are left-invariant vector fields  on the Heisenberg group
satisfying the commutation relations  $[V_j,V_{j+n}]=U$ for $1\le j\le
n$, and where $A=(a_{jk})$ is a complex symmetric matrix with
semi-definite real part. Our results widely  extend all of the results
for the case of non-real, semi-definite matrices
$A$  known to date, in  particular those  obtained recently jointly with
F.~Ricci under Sj\"ostrand's cone  condition. They are achieved by
showing that an  integration by parts argument, which had been applied in
different forms in previous articles, ultimately allows for a reduction to
the case of operators $L$ whose associated Hamiltonian has a purely real
spectrum.  Various examples are  given in order to indicate the potential
scope of this approach and to illuminate some further conditions that
will be introduced in the article.

\end{abstract}

\section{Introduction}\label{introduction}

Consider the standard basis of left-invariant vector 
fields on the Heisenberg group
$\HH_n$, with coordinates $(x,y,u)\in\RR^n\times\RR^n\times \RR$:

\begin{eqnarray} \label{1a}
&&X_j=\de_{x_j}-\frac 12 y_j\de_u\ ,\qquad j=1,\cdots,n\ ,\nonumber\\ 
&&Y_j=\de_{y_j}+\frac 12 x_j\de_u\ ,\qquad j=1,\cdots,n\ ,\\
&&U=\de_u .\nonumber
\end{eqnarray}

We also write $V_1,\dots,V_{2n}$ for
$X_1,\dots,X_n,Y_1,\dots,Y_n$ (in this order), and, consistently,
$v=(x,y)\in\RR^{2n}$.

Given a $2n\times 2n$ complex symmetric matrix
$A=(a_{jk})$, set
\begin{equation}\label{1b}
\L_A=\sum_{j,k=1}^{2n}a_{jk}V_jV_k\ ,
\end{equation}
and, for $\alpha\in\CC$,
\begin{eqnarray}\label{1c}
\L_{A,\alpha}=\L_A+i\alpha U\ .
\end{eqnarray}
These operators can be characterized as the
second order left-invariant differential  operators on $\HH_n$ which
are homogeneous of degree 2 under the automorphic dilations 
$(v,u)\mapsto (\delta v,\delta^2 u)$.

It is the goal of this article to devise large classes of operators 
$\L_{A,\alpha}$ with non-real coefficient matrices $A$ that are locally solvable, extending and unifying 
in this way results from previous articles, including 
\cite{demari-peloso-ricci}, \cite{mueller-peloso-ricci-1},
\cite{mueller-peloso-ricci-2} and \cite{mueller-ricci-cone}. The 
case of real coefficients had been treated in a complete way in 
\cite{mueller-ricci-2}, see also \cite{mueller-ricci-2step}. 
 For background information on the problem of local solvability for 
such operators and further references,  see \cite{mueller-ricci-survey},
\cite{mueller-ricci-cone}. 

In order to motivate the conditions on the coefficient matrix $A$ 
that we shall impose, let me point out that there are various indications 
that,  at least in ''generic'' situations of sufficiently 
high dimension, 
local solvability of $\L_{A,\alpha}$ will only occur if the principal 
symbol $p_A$ satisfies a sign condition, i.e if there exists some 
$\theta\in\RR$ such that $\Re (e^{i\theta} p_A)\ge 0$ (see e.g. 
\cite{mueller-peloso-ricci-2},
\cite{mueller-peloso-nonsolv} and the examples to follow). We know of
locally  solvable operators $\L_{A,\al}$  which do not satisfy a 
sign condition (see 
 \cite{mueller-karadzhov}, \cite{mueller-peloso-nonsolv}), but these examples effectively only occur in 
low dimensions.
For $\theta=0$ (which we may then assume without loss of generality), 
 the sign condition means that $\L_{A,\alpha}$ is a dissipative differential 
operator, or, equivalently, that $\Re A\ge 0.$ The latter condition is what 
we shall assume throughout the paper. This condition is considerably 
weaker than Sj\"ostrand's cone condition (see \cite{sjoestrand}), which 
was imposed in \cite{mueller-ricci-cone} and which for the operators 
\eqref{1c} just means that there is  a constant $C>0$ such that 
$|\Im A|\le C\,\Re A.$

The operators $\L_{A,\al}$ have double characteristics, and for such 
operators it is known that it is not only the principal symbol that governs 
local solvability, but that also  the subprincipal symbol in 
combination with the Hamiltonian mappings associated with doubly characteristic 
points plays an important role. Due to the translation invariance of 
our operators and the symplectic structure that is inherent in the 
Heisenberg group law, these Hamiltonians are essentially encoded in the 
Hamiltonian $S\in \sp(n, \CC),$ associated to the coefficient matrix 
$A$ by the relation $S:=-AJ$ (see e.g. \cite{mueller-ricci-survey}). 
Here, $J$ denotes the matrix $J:=\begin{pmatrix} 0 & I_n\cr
	-I_n & 0\end{pmatrix}.$ In order to emphasize the central role played 
by $S$, we shall therefore also denote $\L_{A,\al}$ by $L_{S,\al}.$ 

One of our main results, Theorem \ref{2A}, states that, under a 
further, natural condition, the question of local solvability of 
the operators  $L_{S,\al}$ can essentially be reduced to the case 
where the Hamiltonian $S$ has only real eigenvalues.

This is achieved by showing that an  integration by parts technic, 
which had been introduced  by R.~Beals and 
P.C.~Greiner in \cite{beals-greiner} and since then been applied in 
modified ways  in  various subsequent articles, e.g. in \cite{mueller-ricci-cone},
 when viewed in the right way,  ultimately allows to 
show that $L_{S,\al}$ is locally solvable, provided that $L_{S_r,\beta}$
 is locally solvable for particular values of $\beta.$ Here, $S_r$ is the 
"part" of $S$ comprising all Jordan blocks associated with real 
eigenvalues.

In combination with some partial results for the case of real eigenvalues,
 this theorem allows to widely extend all the results known to date 
for operators $\L_{A,\al}$ with non-real coefficient matrices $A$  
(see Theorems \ref{2E},\ \ref{2G}).   Moreover, 
we believe that our proofs have become simplier compared e.g. to 
\cite{mueller-ricci-cone}, because of the new structural insights given by 
Theorem \ref{2A}. 

\bigskip
The article is organized as follows. Section \ref{results} contains the 
basic notation, which is mostly taken from \cite{mueller-ricci-cone}, and the statements 
of the main results. Moreover, we present various examples which help 
to illustrate some additional conditions imposed in our theorems.

The preparatory material on the algebraic properties of symplectic 
transformations needed in the proofs of our main results is collected 
in Section \ref{algebra}. 

Section \ref{examples} is devoted to the discussion of examples, including 
those mentioned in Section~\ref{results}.  

In Section \ref{gaussians}, following some line of thoughts in
 \cite{beals-greiner}, we derive explicit formulas for the 
one-parameter semigroups $\{e^{t\L_{A,\alpha}}\}_{t>0}.$ Since 
$\L_{A,\alpha}$ is assumed to be dissipative, these semigroups 
 do exist, at least as linear contractions
on $L^2(\HH_n).$  For the case where $\Re A$ is strictly positive definite,
 such formulas had been established in \cite{mueller-ricci-cone} by 
means of the oscillator semigroup, introduced by R.~Howe in \cite{howe}.
 Starting from formula (5.10) in Theorem 5.2, \cite{mueller-ricci-cone},
 we extend its range of validity to arbitrary matrices $A$ with 
$\Re A\ge 0,$ by adapting some limit arguments from \cite{hoermander-mehler} to 
the setting of twisted convolution operators. We should like to mention 
that the main result in this section, Theorem \ref{5E}, could also be 
derived directly from Theorem 4.3 in \cite{hoermander-mehler} by means 
of the Weyl transform, which relates twisted convolution operators to 
pseudo-differential operators in the Weyl calculus. We prefer to 
present our approach, nevertheless, since we believe that the approach 
through the oscillator semigroup is somewhat simpler than in 
\cite{hoermander-mehler}  and \cite{beals-greiner}.

We also study the  analytic extension of our formulas to the case 
of arbitrary complex matrices $A$ and complex time parameter $t.$ 
This will be useful in the proof of  Theorem \ref {2A}, which will 
be given Section \ref{reduction}, in that it allows to use complex symplectic 
changes of coordinates in some situations.

Finally, in Section \ref{real}, we prove some partial results 
on the case  where $S$ has purely real spectrum, which in combination 
with results from \cite{mueller-ricci-2} lead to Theorem \ref{2E}. 
Moreover, in Proposition \ref{7C}, we shall prove by means of a representation 
theoretic criterion from \cite{mueller-peloso-nonsolv} that the operators 
from Example \ref{2F} are always locally solvable. This result shows that 
the  simple minded approach used in Proposition \ref{7A}  is rather limited,
and that new ideas will be need in order to obtain a better understanding 
of  the case where $S$ has  purely real spectrum.

\section{Statement of the main results}\label{results}

In order to emphasize the symplectic structure on $\RR^{2n}$ which is 
implicit in (1.1), and at the same
time to provide a coordinate-free approach, we shall adopt the notation from
\cite{mueller-ricci-cone} and work within the setting of an arbitrary
$2n$-dimensional real vector space $V$, endowed with a symplectic form $\sigma$.
The extension of
$\sigma$ to a complex symplectic form on $V^\CC$, the complexification of $V$, will
also be denoted by
$\sigma$.

If $Q$ is a complex-valued symmetric form on $V$, we shall often view it as a
symmetric bilinear form on $V^\CC$, and shall denote by $Q(v)$ the quadratic form
$Q(v,v)$. $Q$ and $\sigma $ determine a linear endomorphism $S$ of $V^\CC$ by
imposing that 
\[
\sigma(v,S w)=Q(v,w).
\]
Then, $S\in \sp(V^\CC,\sigma)$, i.e.
\begin{equation}\label{2a}
\sigma (S v,w)+\sigma(v,Sw)=0.
\end{equation}
$S$ is called the {\it Hamilton map} of $V$.
Clearly, $S$ is real, i.e. $S\in \sp(V,\sigma)$, if $Q$ is real.

If $T:U\to W$ is a linear homomorphism  of real or complex vector spaces, we 
shall denote by $\t T:W^*\to U^*$ the transposed homomorphisms between the dual
spaces $W^*$ and $U^*$ of $W$ and $U$, respectively, i.e.
\[
(\t T w^*)(u)=w^*(Tu),\ u\in U, w^*\in W^*.
\]
As usually, we shall identify the bi-dual $W^{**}$ with $W$. 

If $Q$ is any bilinear form on $V$ (respectively $V^\CC$), there is a unique linear 
map $\Q :V\to V^*$ (respectively $\Q:V^\CC\to (V^\CC)^*$) such that 
\begin{equation}\label{2b}
(\Q v)(w)=Q(v,w),
\end{equation}
and $\Q$ is a linear isomorphism if and only if $Q$ is non-degenerate. In
particular,
 the map $\J:V^\CC\to (V^\CC)^*$, given by 
\begin{equation}\label{2c}
(\J v)(w)=\sigma(w,v),
\end{equation}
is a linear isomorphism, which restricts to a linear isomorphisms from $V$ to
 $V^*$, also denoted by $\J$. We have $\t \J=-\J$, so that (\ref{2a}) can be read
as 
\begin{equation}\label{2d}
\J S+\t S\J=0.
\end{equation}
Moreover, if the form $Q$ in (\ref{2b}) is symmetric, then $\t \Q=\Q$ and $\Q=\J S$.

By composition with $\J$, bilinear forms on $V$ can be transported to $V^*$, e.g. 
we put
\[
\sigma^*(\J v,\J w):=\sigma(v,w),\quad  Q^*(\J v,\J w):=Q(v,w).
\]
In analogy with (\ref{2b}) and (\ref{2c}), we obtain maps from $V^*$ to $V$ 
(respectively from $(V^\CC)^*$ to $V^\CC$) which satisfy the following identities:
\begin{equation}\label{2e}
\J^*=-\J^{-1},\  \Q=\t \J \Q^* \J=-\J \Q^* \J,\ S^*=\J S \J^{-1}=-\t S.
\end{equation}

The canonical model of a symplectic vector space is $\RR^{2n}$, with symplectic
form $\sigma (v,w)=\t v Jw=:\l v,w\r$, where 
\[
J:=\begin{pmatrix} 0 & I_n\cr
	-I_n & 0\end{pmatrix}.
\]
Identifying also the dual space with $\RR^{2n}$ (via the canonical inner product on 
$\RR^{2n}$), we have $\J v=Jv$. Moreover, of course $(\RR^{2n})^\CC=\CC^{2n}$.

If a general symmetric form $Q$ is given by $Q(v,w)=\t v Aw$, where $A$ is a symmetric matrix,
we have the following formulas:
\begin{equation}\label{2f}
\Q v=Av,\ Sv=-JAv,\ S^* v=-AJv.
\end{equation}
These formulas apply whenever we introduce coordinates on $V$ adapted to a symplectic basis
 of $V$, i.e. to a basis $X_1,\dots,X_n,Y_1,\dots,Y_n$ such that 
\[
\sigma(X_j,X_k)=\sigma(Y_j,Y_k)=0,\ \sigma(X_j,Y_k)=\delta_{jk}
\]
for every $j,k$. Observe that in $V^*$, the dual of a symplectic basis is
symplectic with respect to $\sigma ^*.$ The Heisenberg group $\HH_V$ built on $V$ is
$V\times \RR$, endowed with the product 
\[
(v,u)(v',u'):= (v+v',u+u' +\frac 1 2 \sigma(v,v').
\]
Its Lie algebra $\h_V$ is generated by the left-invariant vector fields
\[
X_v=\partial_v+\frac 1 2 \sigma(\cdot,v)\partial_u,\quad v\in V.
\]
 The Lie brackets are given by $[X_v,X_w]=\sigma(v,w)U$, with $U:=\partial_u$.

We regard the formal expression (1.2) defining the operator ${\cal L}_A$ as an element 
of the symmetric tensor product $\frak S^2(V^\CC)$ (with $V^\CC=\CC^{2n})$, hence as
a complex symmetric bilinear form $Q^*$ on $(V^\CC)^*$. With this notation, the
Hamilton map $S^*$ of $Q^*$ is 
\begin{equation}\label{2g}
S^* v=-JAv,
\end{equation}
and the Hamilton map of the corresponding form $Q$ on $V^\CC$ is 
\begin{equation}\label{2h}
Sv=-AJv.
\end{equation}
Since the solvability of $\L_{A,\alpha}$ is closely connected to the spectral
properties of the associated Hamilton map, we shall also write
\[
\L_A=:L_S,\quad \L_{A,\al}=:L_{S,\al}.
\]
We remark that $[L_{S_1},L_{S_2}]=-2\ L_{[S_1, S_2]}U$.
The following structure theory for elements $S\in \sp(V^\CC,\sigma)$ will be
important (see e.g. \cite{hoermander-mehler}, \cite{mueller-ricci-cone}).

If $\spec S\subset \CC$ denotes the spectrum of $S$, then $-\lambda \in \spec S$ whenever $\lambda \in \spec S$. 
Moreover, if $V_\la$ denotes the generalized eigenspace of $S$ belonging to the
eigenvalue $\lambda$, then
\begin{equation}\label{2i}
\sigma(V_\la,V_\mu)=0,\quad \mbox{ if } \la +\mu \ne 0.
\end{equation}
In particular, $V_\la$ and $V_{-\la}$ are isotropic subspaces with respect to the symplectic form $\sigma$, and $V_\la \oplus V_{-\la}$ is a symplectic subspace of $V^\CC$, if $\la \ne 0$, while $V_0$ is symplectic too. We 
thus obtain a decomposition of $V^\CC$ as a direct sum of symplectic subspaces
which are $\sigma$-orthogonal:
\begin{equation}\label{2j}
V^\CC=V_0\oplus \sum_{\la \ne 0}^{\quad\oplus} V_\la \oplus V_{-\la}.
\end{equation}
Here, the summation takes place over a suitable subset of $\spec S$. Notice that the
decomposition above is also orthogonal with respect to the symmetric form
$Q(v,w)=\sigma(v,Sw)$,  since the spaces $V_\la$ are $S$-invariant. (\ref{2j})
induces an orthogonal decomposition
\begin{equation}\label{2k} 
V^\CC=V_r\oplus V_i,
\end{equation}
where $V_r:=\sum\limits_{\la \in \RR\cap \spec S}^{\quad \oplus}V_\la$ and 
$V_i:=\sum\limits_{\mu\in (\CC\setminus \RR)\cap \spec S} ^{\quad\oplus } V_\mu$.
Correspondingly, $S$ decomposes as 
\begin{equation}\label{2l} 
S=S_r+S_i,
\end{equation}
where we have put $S_r(u+w):=S(u)$, $S_i(u+w):=S(w)$, if $u\in V_r$ and $w\in V_i$. 
Then also $S_r$ and $S_i$ are in $\sp(V^\CC,\sigma)$, and $S_r$ respectively $S_i$
corresponds to the Jordan blocks of $S$ associated with real eigenvalues
respectively non-real eigenvalues.

Next, $S$ can be uniquely decomposed as $S=D+N$ such that $D$ is semisimple, $N$ is nilpotent and $DN=ND$. The endomorphisms $D$ and $N$  in this Jordan decomposition of $S$ 
are polynomials in $S$, and it is known from general Lie theory that $D,N\in
\sp(V^\CC,\sigma)$ (see e.g. \cite{borel}).

 This can also be seen directly. For this
purpose, we may assume without loss of generality  that $V^\CC=V_\la \oplus
V_\la,$ for some
$\la\ne 0$ (the case $V^\CC=V_0$ is obvious). Then, if $v_\la,w_\la $ are in
$V_\la$ and
$v_{-\la},w_{-\la}$ in $V_{-\la}$, we have
\begin{eqnarray*}
\lefteqn{\sigma (D(v_\la +v_{-\la}),w_\la+w_{-\la})=\sigma(\la v_\la - 
\la v_{-\la},w_\la+w_{-\la})}\\
&&=\sigma(v_\la,\la w_{-\la})+\sigma(v_{-\la},-\la w_\la)\\
&&=-\sigma(v_\la +v_{-\la},D(w_\la +w_{-\la})).
\end{eqnarray*}
Applying the Jordan decomposition to $S_r$, we  can uniquely write
\begin{equation}\label{2m}
S_r=D_r+N_r,
\end{equation}
with $D_r$ semisimple, $N_r$ nilpotent and $D_rN_r=N_rD_r$.

Our first main result is a reduction theorem, which allows in many cases to reduce the 
question of local solvability of $L_{S,\al}$ to essentially  the same question for
the operator $L_{S_r,\beta}$, for particular values of $\beta\in \CC$. Its proof is
based on an integration by parts argument, variants of which had already been used
in \cite{beals-greiner} as well as in several subsequent articles, e.g. in
\cite{mueller-peloso-ricci-1}, \cite{mueller-ricci-cone}. We believe that our
approach reveals more clearly and conceptually the potential range of validity of
such technics of integration by parts, by showing that they allow a reduction to the
study of the operators
$L_{S_r,\beta}$.

A main obstruction to applying this technic is the fact that the spaces 
$V_\la\oplus V_{-\la}$, for $\la \in \RR\setminus \{0\}$, are in general not
invariant under complex conjugation. This had already been observed by
L.~H\"ormander \cite{hoermander-mehler}, who gave an example for the related case
$\la=0$, and we shall give further examples in Section \ref{examples}.

We shall therefore mostly work under the following hypothesis:
\begin{equation}\tag{R}\label{R}
V_\la\oplus  V_{-\la}\text{ is invariant  under complex  conjugation, for every } \la \in \RR\setminus \{ 0\}.
\end{equation}

We write 
\[
\spec S_i=\{\pm \om_1,\dots,\pm \om_{n_1}\}\subset \CC\setminus \RR,
\]
where the eigenvalues $\pm  \om_j$ are listed with their multiplicities, and where 
\begin{equation}\label{2n}
\nu_j:=\Im \om_j>0, \ j=1,\dots,n_1.
\end{equation}
We also put 
\[
\nu:=\sum_{j=1}^{n_1} \nu_j,\quad \nu_{\min} :=\min_{j=1,\dots,n_1} \nu_j\ .
\]

\begin{theorem}\label{2A}
Assume that $\Re Q_S\ge 0$, that $S_i\ne 0$ (i.e. $\nu>0$), and that condition
(\ref{R}) is satisfied. Then, the following holds:
\be
\item[(i)]$L_S +i\al U$ is locally solvable (and even admits a tempered 
fundamental solution), if $|\Re\al|<\nu$.
\item [(ii)] If $M>0$, and if $|\Re \al|<M$, then
$L_S+i\al U$ is locally solvable, provided that  $L_{S_r} +i(\al \pm
\sum\limits_{j=1}^{n_1}  (2k_j+1)i\om_j)U$ is locally solvable, for every
$n_1$-tupel
$(k_1,\dots,k_{n_1})
\in \NN^{n_1}$ such that 
$\sum\limits_{j=1}^{n_1} (2k_j+1)\nu_j< M$.
\ee
\end{theorem}

In view of Theorem \ref{2A}, it thus becomes a major task to understand local solvability when $S$ has purely real spectrum.

There are various indications that, at least in sufficiently high dimensions,
 $L_{S,\al}$ may not be locally solvable, unless $S$ satisfies the following sign
condition:
\begin{equation}\label{2o}
\Re(e^{i\theta} Q_S)\ge 0 \ \mbox{ for some } \theta \in \RR
\end{equation}
(see e.g. \cite{mueller-peloso-nonsolv}, and also the examples to follow). We shall
therefore assume that \eqref{2o} holds for $S_r$, and  even that, without loss of
generality, 
$\Re Q_{S_r}\ge 0$. The following proposition gives a sufficient condition for this
to hold.

\begin{proposition}\label{2B}
Assume that property (\ref{R}) is satisfied and that $N_r^2=0$, where $N_r$ denotes the nilpotent part 
in the Jordan decomposition of $S_r$. Then $\Re Q_S\ge 0$ implies $\Re Q_{S_r}\ge
0$.
\end{proposition}

Consider the following examples, which shed some more light on the conditions in
Proposition \ref{2b}.

\begin{example}\label{2C}

{\rm On $\HH_3$, consider
\[
L_S:=Y_1^2+X_3^2+2X_3Y_1+Y_3^2+2i(X_1Y_2-X_2Y_1-Y_2Y_3).
\]
It is clear that $\Re Q_S\ge 0$. However, we will show in Section \ref{examples}
that neither $\Re Q_{S_r}\ge 0$, nor $\Re Q_{S_i}\ge 0$, even though $N_r$ is
2-step nilpotent. Even worse, by H\"ormander's criterion (H), one checks that
$L_{S_r}+F$ und also
$L_{S_i}+F$ are not locally solvable, for every first order differential operator
$F$ with smooth coefficients. This shows that Proposition \ref{2B} fails to be true
without property (\ref{R}), and that the entire approach  in Theorem \ref{2A} will
in general break down, if (\ref{R}) is not satisfied. Theorem \ref{2A} (i) remains 
nevertheless valid in this example, i.e. $L_S+i\al U$ is locally solvable for 
$|\Re \al|<1$ (see Remark \ref{rem6} following Proposition \ref{6A}).We do not
know what  happens if $|\Re\al|\ge 1.$}
\end{example}

\begin{example}\label{2D}
{\rm On $\HH_3$, consider
\[
L_S:=X_2^2+X_3^2+Y_3^2+2i(X_1Y_2+bX_2Y_3),\quad b\in \RR\setminus \{0\}.
\]

Again, we have $\Re Q_S\ge 0$. We shall see that $S_r=N_r$ is nilpotent of step 4
in this example, and that $\Re Q_{S_r}\ge 0$ is not satisfied (not even \eqref{2o}).
Again, H\"ormander's condition (H) is satisfied by $L_{S_r}$, so that $L_{S_r} +F$
is not  locally solvable, for every first order term $F$. Notice that in this
example property (R) does hold, since $S_r=N_r$, which shows that the conclusion
of Proposition \ref{2b} will in general not hold, if $N_r$ is nilpotent of step
higher than 2, even under property (R).

In our study of $L_{S_r}$, we shall therefore restrict ourselves to the case where
$N_r^2=0$.

Let us thus assume, for a moment, that $\spec S\subset \RR$, i.e. $S=S_r$, and
that $N^2=0$,
 where $S=D+N$ is the Jordan decomposition of $S$. We write $D=D_1+iD_2$, where
$D_1=\Re D$, $D_2=\Im D$. If we assume that $D$ satisfies the following hypothesis

\begin{equation}\tag{C}\label{C}
[D_1,D_2]=0,
\end{equation}
then we can discuss local solvability of $L_S+i\al U$ in a complete way, even without 
assuming (R),
 by means of certain a priori estimates and the results in \cite{mueller-ricci-2}.}
\end{example}

\begin{theorem}\label{2E}
Assume that $S\ne 0$ has purely real spectrum, $\Re Q_S\ge 0$, $N^2=0$, and that
property (C) is satisfied. Then the following holds true:
\be
\item[(i)] If $\Re S\ne 0$ or $\Re\al \ne 0$, then $L_S+i\al U$ is locally solvable.
\item[(ii)] If $\Re S=0$, $\Re\al =0$ and $N\ne 0$, then $L_S+i\al U$ is locally
solvable.
\item[(iii)]If $\Re S=0,\ \Re\al =0$ and $N=0$ then,
after applying a suitable automorphism of $\HH_n$ leaving the center fixed, $L_S$ takes the form 
\begin{equation}\label{2p}
L_S=\sum_{j=1}^n i\la _j(X_j^2+Y_j^2),
\end{equation}
with $\la _1,\dots,\la_n\in \RR.$ In particular, $\spec S=\{\pm \la_1,\dots,\pm
\la_n\}$. Then $L_S+i\al U$ is locally solvable if and only if there are constants
$C>0$ and $M\in \NN$, such that 
\begin{equation}\label{2q}
|\al\pm \sum_{j=1}^n (2k_j+1)i\la_j|\ge C(1+|k|)^{-M},
\end{equation}
for every $k=(k_1,\dots,k_n)\in \NN^n$.
\ee
\end{theorem}

Consider the following example.

\begin{example}\label{2F}
{\rm On $\HH_2$, let
\[
L_S:=(m+c_1) Y_1^2 +(m-c_1)Y_2^2+2c_2Y_1Y_2+2i(X_1Y_2-X_2Y_1),
\]
where we assume that $c_1^2+c_2^2\ne 0$. It is easy to see that for
$m\ge \sqrt{c_1^2+c_2^2}$ one has $\Re Q_S\ge 0$. 
Moreover, we shall show that $D$ and $N$ are given by the block matrices
\[
D=\begin{pmatrix}
iJ & 0\\
C & iJ
\end{pmatrix},
\quad N=\begin{pmatrix}
0 & 0\\mI & 0\end{pmatrix},
\]
where $J=\begin{pmatrix}0 & 1\\-1 & 0\end{pmatrix}, 
C=\begin{pmatrix}c_1 & c_2\\c_2 & -c_1\end{pmatrix}$ and $I=\begin{pmatrix}1 & 0\\0
& 1\end{pmatrix}$, and that $D^2=I$, so that $N^2=0$ and $\spec S=\{-1,1\}$. One
checks immediately that property (C) is not satisfied in this example.
Nevertheless, we shall prove in Section \ref{real}, Proposition \ref{7C}, by means
of some explicit computations based on the group Fourier transform on $\HH_2$, that
$L_S+i\al U$ is locally solvable, for every $\al \in \CC$.

This example indicates that condition (\ref{C}) may not be necessary in Theorem
\ref{2E}.
 We shall further comment on Example \ref{2F} in Section \ref{examples}.

\medskip
The conditions in Theorem \ref{2E} are of course  rather restrictive, but at
present we do not know of any approach which would allow to discuss much wider 
classes of operators $L_S+i\al U$, with $\al \in \CC$ and $\spec S\subset \RR,\  \Re
Q_S\ge 0$, even if $S$ is nilpotent.  A first, useful step towards a better
understanding of such operators might be a classification of normal forms of
matrices $S$ satisfying $S^2=0$, along the lines of \cite{mueller-thiele-forms}.

Nevertheless, Theorem 2.5 in combination with Theorem 2.1 immediately gives the
subsequent theorem. It contains and widely extends, in combination with Theorem
\ref{2E}, all of the  positive results on local solvability
which have been obtained hitherto under  the sign condition
\eqref{2o} in the ``non-real'' case (but for the discussion of the exceptional
values arizing in (ii)).}
\end{example}

\begin{theorem}\label{2G}
Let $S\in \sp(n,\CC)$ be such that $\Re Q_S\ge 0$. Assume further that $S$ has at
least one non-real eigenvalue and satisfies property (\ref{R}), and that $N_r^2=0$
and $\Re D_r=0$, where $S_r=D_r+N_r$ is the Jordan decomposition of $S_r$. Let
$\om_1,\dots,\om_{n_1}$ be as in Theorem~\ref{2A}. Then the following holds:
\be
\item [(i)] $L_S+i\al U$ is locally solvable for every $\al \in \CC$, provided $S_r\ne 0$.
\item[(ii)] If $S_r=0$, then $L_S+i\al U$ is locally solvable for all values of 
$\al$ in $\CC\setminus \E_S$, where $\E_S$ is the following set of exceptional
values:
\[
\E_S:=\{\pm \sum_{j=1}^{n_1} (2k_j+1)i\om _j: k_1,\dots, k_{n_1}\in \NN\}.
\]
\ee
\end{theorem}

\begin{remarks}\label{2H}
{\rm 
\noi (a) Assume that $L_S$ satisfies the cone condition in the sense of Sj\"ostrand
and H\"ormander, i.e. there exists a constant $C>0$, such that 
\begin{equation}\label{2r}
|\Im Q_S(v)|\le C|\Re Q_S (v)|\quad \forall v \in \RR^{2n}.
\end{equation}
Then it is known that $S$ has at most one real eigenvalue, namely 0, that $S_r=N_r$
is 2-step nilpotent, and that $S_r=0$ if and only if $\ker S$ is symplectic (see
\cite{mueller-ricci-cone}, Lemma 3.6).

Thus, Theorem \ref{2G} contains Theorem 2.2 in \cite{mueller-ricci-cone}, except for
the proof of non-solvability for the exceptional values $\al\in \E_S$, in case that
$S_r=0$.

\noi (b) In the presence of property (\ref{R}), the seemingly stronger condition
$\Re D_r=0$ is in fact equivalent to the condition (\ref{C}) for $S_r$, i.e. $[\Re
D_r,\Im D_r]=0$ (see Corollary \ref{3H}). }
\end{remarks}

\setcounter{equation}{0}
\section{On the algebraic structure of $S$}\label{algebra}

Let $S\in \sp(V^\CC,\sigma)$ be the Hamilton map associated to the
quadratic form $Q$ on $V$, which consequently  we sometimes will also
denote by $Q_S$. Our general assumption is that
\begin{equation}\label{3a}
\Re Q \ge 0 \text{ on } V.
\end{equation}
This condition is equivalent to the following condition on the $\CC$-linear
 extension of $Q$ to $V^\CC\times V^\CC:$
\begin{equation}\label{3b}
\Re Q(z,\overline z)\ge 0\quad \forall z\in V^\CC.
\end{equation}
Here, complex conjugation in $V^\CC$ is meant with respect to the real
form $V$ of $V^\CC,$ i.e. for $z=v+iw\in V^\CC$ we set $\overline
z:=v-iw$.

In the sequel, we shall often indicate the real part of a linear map or form by an index 1, the imaginary part by an index~2. For instance,
\[
S=S_1+iS_2,\quad Q=Q_1+iQ_2.
\]
The following proposition is due to H\"ormander (\cite{hoermander-mehler},
Proposition 4.4).

\begin{proposition}\label{3A}
Assume that $Q_1=\Re Q\ge 0$.
Then the following hold true for every real eigenvalue $\la$ of $S$:
\begin{equation}\label{3c}
\overline{\Ker (S-\la)} =\Ker (S+\la);
\end{equation}

\begin{equation}\label{3d}
S_1\Ker (S\pm \la)=0.
\end{equation}

In particular, $\Ker  (S-\la)\oplus \Ker (S+\la)$, $0\ne \la \in \RR$, 
is the complexification of its intersection with $V$, and so is $\Ker S$.
\end{proposition}

In \cite{hoermander-mehler} H\"ormander gives an example which
shows that the space $V_\la+V_{-\la},\ \la \in \RR$, is in general not
self-conjugate, namely for the case $\la=0$. We remark that our Example
2.3 presents a corresponding example for the case $\la =1$. The following
result is obvious.

\begin{lemma}\label{3B}
Assume that $\la\in \RR\setminus {0}$, and that $V_\la\oplus V_{-\la}$ is 
self-conjugate. Then $V_\la^\RR:=V\cap (V_\la\oplus V_{-\la})$ is an
$S_1$ and $S_2$-invariant real subspace of $V$ such that $V_\la\oplus
V_{-\la}=V_\la^\RR+iV_\la^\RR$, and $\Re Q\ge 0$ on $V_\la^\RR$, if
\eqref{3A} holds.  In particular, $\Re Q\ge 0$ on $\sum\limits_{\la\in
\RR\setminus \{0\}} V_\la^\RR$, if $S$ satisfies condition \eqref{R}.
\end{lemma}

\begin{lemma}\label{3C}
Assume that condition \eqref{3a} holds. Then $\Re Q_{S_r}\ge 0$ if and
only if
$\overline{S(V_r)}\subset V_r$.
\end{lemma}

\noindent{\bf Proof.} We have $\Re Q_{S_1}\ge 0$ if and only if
\[
\Re \sigma (\overline v,S_r v)\ge 0\quad \forall v\in V^\CC.
\]
Let $v=u+w, u\in V_r,w\in V_i$. Then 
\begin{eqnarray*}
\Re \sigma (\overline v,S_rv)&=&\Re \sigma (\overline u+\overline w,Su)
= \Re Q_S(u,\overline u)+\Re \sigma (\overline w, Su).
\end{eqnarray*}

This is non-negative for every $v\in V^\CC$ if and only if $\Re \sigma (\overline w,Su)=0$ for every $u\in V_r,w\in V_i$. Since $V_r$ and $V_i$ are complex vector spaces, this means that $\sigma(\overline w,z)=0$ or, equivalently, $\sigma(w,\overline z)=0$ for every $w\in V_i$ and $z\in S(V_r)$. $V_i$ being the orthogonal complement of $V_r$ (w.r. to $\sigma$), the latter condition is equivalent to $\overline {S(V_r)}\subset V_r$.

\hfill Q.E.D.
\bigskip

We can now easily prove Proposition \ref{2B}.
Since $V_r=V_0\oplus \sum\limits_{\la\in \RR\setminus \{0\}} (V_\la \oplus V_{-\la})$, where each of the occuring subspaces is $S$-invariant, and where $V_\la \oplus V_{-\la}$ is self-conjugate for $\la\ne 0$ (because of property (\ref{R})),
we only have to show that $\overline{S(V_0)}\subset V_0$, in order to
apply Lemma \ref{3C}. But, $S(V_0)=N_r(V_0)\subset \Ker N_r \cap V_r=\Ker
S$, and
$\overline {\Ker S}=\Ker S$ by Proposition \ref{3A}, so that
$\overline{S(V_0)}
\subset \Ker S\subset V_0$, which completes the proof.

\hfill Q.E.D.
\bigskip

The form $\Re Q_{N_r}$ is always semi-definite, if $N_r^2=0$, as the following result shows.

\begin{lemma}\label{3E}
If $N_r^2=0$, then (\ref{3a}) implies $\Re Q_{N_r}\ge 0$ (even without property (\ref{R})).
\end{lemma}

\noindent{\bf Proof.} Decompose $v\in V^\CC$ as $v=w+\sum\limits_{\la \in \RR\cap \spec S} v_\la$, where $w\in V_i$ and $v_\la \in V_\la$. Then 
$N_rv=\sum\limits_{\la} (S-\la)v_\la$, where $(S-\la)v_\la\in \Ker
(S-\la)$, since $N_r^2=0$. By Proposition \ref{3A}, we have then 
$\overline {(S-\la)v_\la}\in \Ker (S+\la)\subset V_{-\la}$, so that 
\begin{eqnarray*}
\lefteqn{\sigma(v,\overline{N_rv})= \sigma(w+\sum_\mu v_\mu,\sum_\la 
\overline{(S-\la)v_\la})}\\
&&=\sum_{\la,\mu} \sigma (v_\mu,\overline{(S-\la)}v_\la)=\sum_\la \sigma(v_\la,\overline{(S-\la)v_\la}),
\end{eqnarray*}
since $\sigma(V_r,V_i)=0$. Moreover,
$\Re \sigma(v_\la,\overline{(S-\la)v_\la})=\Re \sigma (v_\la,\overline{Sv_\la})$, since 
${\Re \sigma(z,\overline z)=0}\quad \forall z\in V^\CC$. We thus obtain 
$$\Re \sigma(v_\la,\overline{(S-\la)v_\la})=\Re \sigma 
(\overline v_\la, Sv_\la)
=\Re Q_S(v_\la,\overline v_\la)\ge 0,
$$
and then
$$\Re Q_{N_r}(v,\overline v)=\Re Q_{N_r}(\overline v,v)=\Re
\sigma(\overline v,N_rv)=\Re (v,\overline {N_rv})\ge 0.
$$

\hfill Q.E.D.
\bigskip

\noindent From now on, we assume that $\spec S\subset \RR$, i.e. $S=S_r$ and $N=N_r$, and that $N^2=0$.
By Lemma \ref{3E}, this implies $\Re Q_N\ge 0$.
We put
\[
W:=N_1(V)+N_2(V)=\{\Re (Nz): z\in V^\CC\},
\]
if $N=N_1+iN_2$, and 
\[
K:=W^\bot =\{v\in V:\sigma(v,w)=0\quad \forall w\in W\}.
\]

\begin{proposition}\label{3F}
If $\spec S\subset \RR, \Re Q_S\ge 0$ and $N^2=0$, then $W$ is an
isotropic subspace of $V$, and $K=\Ker N_1\cap \Ker N_2\supset W$.

Moreover, $S_1(K)=0,\  S_2(K)\subset K, \ S_2(W)\subset W$, so that in
particular also\hfill\newline $D_1(K)=0, \ D_2(K)\subset K$ and
$D_2(W)\subset W$. 

Finally
\begin{equation}\label{3e}
K^\CC=\sum_{\la \in \spec S} \Ker S_\la
\end{equation}
and 
\begin{equation}\label{3f}
W=\sum_{\la \in \spec S,\la\ge 0} W_\la,
\end{equation}
where $W_\la:=W\cap (\Ker (S-\la)+\Ker(S+\la))$.
\end{proposition}

\noindent{\bf Proof.} Since $N^2=(N_1+iN_2)^2=0$, we have 
\begin{equation}\label{3g}
N_1^2=N_2^2 \text{ and } N_1N_2+N_2N_1=0.
\end{equation}
Denote by $V_\mu^1$ the generalized eigenspace of $N_1$ corresponding to $\mu\in \spec N_1\subset \CC$.
(\ref{3g}) implies that 
$N_2(V_\mu^1)\subset V_{-\mu}^1$. Moreover, since $Q_{N_1}\ge 0$, all eigenvalues $\mu$ of $N_1$ are purely imaginary.
Indeed, if $N_1v=\mu v, v\ne 0$ and $\mu\ne 0$, then 
$Q_{N_1}(v,\overline v)\ne 0$
 (since otherwise $N_1v=0$), and $Q_{N_1}(v,\overline v)=\sigma(\overline
v,N_1v)=\mu \sigma(v,\overline v)$, where $\sigma(v,\overline v)\in
i\RR$. This implies $\mu\in i\RR$. Thus, either $\mu=0$, or
$\mu=i\nu,\nu\in
\RR\setminus \{0\}$. We shall exclude the second possibility.

Namely, if $\nu\in \RR\setminus\{0\}$, then $V_{i\nu}^1\oplus V_{-i\nu}^1$
is an $N_1$ and $N_2$-invariant symplectic subspace of $V^\CC$, obviously
also invariant under complex conjugation. Thus $V_{i\nu}^1\oplus
V_{-i\nu}^1$ is the complexification of its intersection with $V$. We may
therefore, for a moment, restrict ourselves to the latter subspace and
assume that
$V^\CC=V_{i\nu}^1\oplus V_{-i\nu}^1$. But, since $Q_{N_1}\ge 0$, we then
have in fact $Q_{N_1}>0$, which means that $Q_N$ satisfies the
cone-condition. Consequently, by \cite{hoermander-mehler}, Lemma 3.2, 
(see also \cite{mueller-ricci-cone}, Lemma 3.1),
$\ker N=0$, a contradiction.

We have shown that $N_1$ has the only eigenvalue 0, and since $Q_{N_1}\ge 0$, by the classification 
of normal forms of quadratic forms on symplectic vector spaces (see e.g.
\cite{hoermander-mehler},Theorem~3.1], we have $N_1^2=0$. By 
(\ref{3g}), then also
$N_2^2=0$.

We will show that indeed 
\begin{equation}\label{3h}
N_1^2=N_2^2=N_1N_2=N_2N_1=0.
\end{equation}
By (\ref{3g}), for $v\in V$ we have 
\begin{eqnarray*}
\lefteqn{Q_{N_1}(N_2v,N_2v)=\sigma(N_2v,N_1N_2v)}\\
&&=-\sigma(N_2,v,N_2N_1v)=\sigma(N_2^2 v,N_1v)=0,
\end{eqnarray*}
which implies that $N_1(N_2 v)=0$, since $Q_{N_1}\ge 0$.
Consequently, $N_1N_2=0$, hence also $N_2N_1=0$.

\eqref{3h} shows that $W=N_1(V)+N_2(V)$ is an isotropic subspace of $V$
which is contained in $\Ker N_1\cap \Ker N_2$. But one sees easily that
$\Ker N_1\cap \Ker N_2=W^\bot =K$.

Let us decompose $V^\CC=\sum_\la ^\oplus V_\la$, where summation is over all $\la\in \spec S\subset \RR$, and let $v=\sum\limits_\la v_\la \in V^\CC$, with $v_\la \in V_\la$. In order to prove (\ref{3e}), we observe that $K^\CC=\Ker N\cap \Ker \overline N$, and that $v\in \Ker N$ if and only if $v_\la\in \Ker (S-\la)$ for every $\la \in \spec S$, i.e. $\Ker N=\sum\limits_\la \Ker (S-\la)$. By (\ref{3c}), 
this space is self-conjugate, so that $\Ker \overline N=\overline{\Ker
N}$, which shows (\ref{3e}). From (\ref{3e}) and (\ref{3d}) we obtain
$S_1(K)=0$.

Next, for $v_\la \in V_\la$, we have $D v_\la=\la v_\la$, hence 
\begin{equation}\label{3i}
N v_\la =(S-\la)v_\la \in \Ker (S-\la).
\end{equation}
Together with (\ref{3c}), this implies 
\[
S(Nv_\la + \overline{Nv_\la})=\la (Nv_\la-\overline{Nv_\la}),
\]
so that $S(\Re (Nv_\la))=i\la \Im (Nv_\la)$. Consequently, $S_2(\Re (Nv_\la))=\la \Im (Nv_\la)\in W$, which shows that $W$ is $S_2$-invariant. Since $S_2\in \sp (V,\sigma)$, this implies that also $K=W^\bot $ is $S_2$-invariant.

Finally, (\ref{3f}) is an immediate consequence of 
 (\ref{3i}) and (\ref{3c}).

\hfill Q.E.D.
\bigskip

\begin{corollary}\label{3G} Under the hypotheses of Proposition \ref{3e},
we have $S_1^2=0$ and $[S_1,S_2]=[D_1,D_2]$.
\end{corollary}

\bigskip

\noindent{\bf Proof.} Since $S_1(K)=0$, i.e. $K\subset \Ker S_1$, we 
have $S_1(V)=(\Ker S_1)^\bot \subset K^\bot=W\subset K$, so that 
\begin{equation}\label{3j}
S_1(K)=0, \ S_1(V)\subset W.
\end{equation}
This implies in particular $S_1^2=0$, and since $N_1(K)=0, \ N_2(V)\subset W$,
 also 
\begin{equation}\label{3k}
D_1(K)=0,\ D_1(V)\subset W.
\end{equation}
Consequently, $D_1^2=0$ and $N_jD_1=D_1N_j=0, \ \text{for} \  j=1,2$.
Moreover,
$[D,N]=0$ implies $[D_1,N_1]=[D_2,N_2]$ and $[D_1,N_2]=-[D_2,N_1]$, so
that in fact $0=[D_1,N_1]=[D_2,N_2]$ and $0=[D_1,N_2]=-[D_2,N_1]$. Since
also $[N_1,N_2]=0$, we find that $[S_1,S_2]=[D_1,D_2]$.

\hfill Q.E.D.
\bigskip

\begin{corollary}\label{3H} Assume that $S$ satisfies the hypotheses of
Proposition \ref{3F} as well as property~(\ref{R}). Then $[D_1,D_2]=0$ if
and only if $D_1=0$.
\end{corollary}

\bigskip

\noindent{\bf Proof.} One implication being trivial, we assume that $[D_1,D_2]=0$.
 Because of \eqref{R}, we can decompose $V$ as $V=\sum\limits_{\la \ne
0}^{\quad\oplus }(V_\la\oplus V_{-\la})\cap V\oplus V_0\cap V$, where all
subspaces in this decomposition are $S_1$-  and $S_2$-invariant,
symplectic and pairwise orthogonal. We may therefore reduce ourselves to
one of these spaces, i.e. we may assume that $\spec S=\{-\la,\la\}$, for
some $\la\in \RR$. The case $\la=0$ being trivial, let $\la \ne 0$. Then,
$D^2=\la^2 I$, hence 
\[
D_1^2-D_2^2=\la^2 I,\ D_1D_2+D_2D_1=0.
\]
But, from (\ref{3k}), we know that $D_1^2=0$, and since $D_1D_2=D_2D_1$, we thus find that $D_1D_2=D_2D_1=0$ and $D_2^2=-\la^2I$. This implies $D_1=0$.

\hfill Q.E.D.
\bigskip

We conclude this section with a result, which shows that the only way that $\Re Q_D$ 
can be semi-definite is that $D_1=0$.

\begin{lemma}\label{3I}
Let $D=D_1+iD_2\in  \sp(V^\CC,\sigma)$ be such that $D$ is semisimple and $\spec D\subset \RR$.
Then $\Re Q_D$ is semi-definite if and only if $D_1=0$.
\end{lemma}

\bigskip

\noindent{\bf Proof.} Let $Q_D=Q_1+iQ_2$. We first observe that if $Q_1$ is semi-definite, i.e. if $Q_1\ge 0$ or $-Q_1\ge 0$, then $Q_1(z,\overline z)=0$ implies that $Q_1(w,z)=0$ for every $w\in V^\CC$, that is, $D_1 z=0$.

Next, if $0\ne v\in V_\la$, we have $Dv=\la v$, hence $Q_1(v,\overline
v)=\la \Re \sigma(v,\overline v)=0$, so that $D_1v=0$. We have thus shown
that $D_1$ vanishes on every eigenspace $V_\la$ of $D$, hence $D_1=0$.

Inversely, if $D_1=0$, then clearly $\Re Q_D=0$ is semi-definite.

\hfill Q.E.D.


\setcounter{equation}{0}
\section{Examples.}\label{examples}

Before we turn to the proofs of our main theorems, we shall discuss some examples,  
including the ones from Section \ref{results}, in order to illustrate the conditions
we  imposed in  our algebraic results of Section \ref{algebra}.

Our first example demonstrates that the conditions in Theorem \ref{2G} are weaker than the cone condition.

\begin{example}\label{4A}
\rm{On $\HH_2$, consider 
\[
L_S:=X_1^2+X_2^2+i(X_2Y_2+Y_2X_2).
\]
Obviously, $L_S$ does not satisfy the cone condition.
But, $S=-AJ$ has one Jordan block 
$\begin{pmatrix}0 & 0\\ 1 & 0\end{pmatrix}$, and one block
$\begin{pmatrix}i &0 \\ 0 & -i\end{pmatrix}$, since $SY_1=-X_1,\ SX_1=0,$
and $SX_2=iX_2,$ $S(Y_2-\frac i 2 X_2)=-i(Y_2-\frac i 2 X_2)$. 
Thus,
$S_r=N_r\ne 0$ and $N_r^2=0$.}
\end{example}

\noi {\bf Example \ref{2C}}
The operator in this example can be written as 
$$
L_S=2i(X_1Y_2-(X_2+iX_3)Y_1)+Y_1^2+Y_2^2
+X_3^2+(Y_3-iY_2)^2.
$$
Observe that 
$\tilde X_1:=X_1,\tilde X_2:=X_2+iX_3,\tilde X_3:=X_3, \tilde Y_1:=Y_1,\tilde Y_2:=Y_2,\tilde Y_3:=Y_3-iY_2$ is a complex symplectic basis of $\CC^6$, and that $L_S$ can be written
\[
L_S=2i(\tilde X_1\tilde Y_2-\tilde X_2\tilde Y_1)+\tilde Y_1^2 +\tilde Y_2^2+\tilde X_3^2 + \tilde Y_3^2.
\]
The matrix $\tilde S$ corresponding to the new basis is thus given by
\[
\tilde S=-\left(
\begin{array}{cc|cc}
	0   & iJ & & \\
	-iJ & I   & \multicolumn{2}{c}{\raise5pt\hbox{{\rm 0}}}\\
\hline
          \multicolumn{2}{c|}{\lower8pt\hbox{{\rm 0}}} & 1 & 0 \\
              && 0 & 1
\end{array}\right)
\cdot
\left(
\begin{array}{cc|cc}
	0   & I& & \\
	-I &  0   & \multicolumn{2}{c}{\raise5pt\hbox{{\rm 0}}}\\
\hline
          \multicolumn{2}{c|}{\lower8pt\hbox{{\rm 0}}} 
	& 0 & 1 \\
              && -1 & 0
\end{array}\right)
=
\left(
\begin{array}{cc|cc}
	iJ   & 0 & & \\
	I & iJ   & \multicolumn{2}{c}{\raise5pt\hbox{{\rm 0}}}\\
\hline
          \multicolumn{2}{c|}{\lower8pt\hbox{{\rm 0}}} & 0 & -1 \\
              && 1 & 0
\end{array}\right),
\]
with respect to the blocks of symplectic coordinates corresponding to
$\tilde X_1,\tilde X_2,\tilde Y_1,\tilde Y_2$ and $\tilde X_3,\tilde Y_3.$
Here,
$J=\begin{pmatrix}0&1\\-1&0\end{pmatrix}$ and
$I=\begin{pmatrix}1&0\\0&1\end{pmatrix}$. The first block of $\tilde S$
has eigenvalues $\pm 1$, the second $\pm i$, hence the first corresponds
to $S_r$ and the second to $S_i$. We thus find that
$L_S=L_{S_r}+L_{S_i}$, with 
\begin{eqnarray*}
L_{S_r}&=& 2i(X_1Y_2-(X_2+iX_3)Y_1)+Y_1^2+Y_2^2\\
&=& Y_1^2+Y_2^2+2X_3Y_1+2i(X_1Y_2-X_2Y_1)
\end{eqnarray*}
and 
\[
L_{S_i}=X_3^2+(Y_3-iY_2)^2=X_3^2+Y_3^2-Y_2^2-2iY_2Y_3.
\]
This shows that we have $\Re Q_S\ge 0$, but neither $\Re Q_{S_r}\ge 0$ nor
$\Re Q_{S_i}\ge 0$, even though $N_r$ is obviously 2-step nilpotent.

In order to see that  $L_{S_r}$ and $L_{S_i}$ satisfy H\"ormander's condition (H), observe that for $S',S''\in \sp(n,\RR)$, 
the Poisson bracket of the principal symbols of $L_{S'}$ and $L_{S''}$ 
corresponds to the principal symbol of the commutator $[L_{S'},L_{S''}]$.

And,
$$[L_{\Re S_r},L_{\Im S_r}]= [Y_1^2+Y_2^2+2X_3Y_1,X_1Y_2-X_2Y_1]
=-2X_3Y_2U,
$$
so that, at the origin, for the operator ($L_R +  \text{first order
term}$), the condition (H) reduces to solving the system 
\[
\eta_1^2+\eta_2^2+2\xi_3\eta_1=0,\quad \xi_1\eta_2-\xi_2\eta_1=0,\quad
\xi_3\eta_2\ne 0.
\]
One solution is given by $\xi_1=\xi_2=\xi_3=1,\ \eta_1=\eta_2=\eta_3=-1$.

Similarly, since $[X_3^2+Y_3^2-Y_2^2,Y_2Y_3]=2X_3Y_2U$, condition (H) for 
the operator $L_{S_i}$ reduces to solving the system
\[
\xi_3^2+\eta_3^2-\eta_2^2=0,\quad \eta_2\eta_3=0,\quad \xi_3\eta _2\ne 0.
\]
A solution is given whenever $\eta_3=0$ and $\xi_3=\eta_2=1$.

\bigskip
\noi {\bf Example \ref{2D}}
For $b\in \RR\setminus \{0\}$, consider on $\HH_3$ the operator 
\[
L_S:=X_2^2+X_3^2 + Y_3^2+2i(X_1Y_2+bX_2Y_3).
\]
Then, clearly $\Re A\ge 0$, and $S=-AJ$ is given by 
\[
S=
\left(
\begin{array}{ccc|ccc}
0&&&0&i&0\\
&1& &0&0&ib\\
&&1&0&0&0\\
\hline
0&0&0&0&&\\
i&0&0&&0&\\
0&ib&0&&&1\end{array}
\right)
\cdot 
\left(
\begin{array}{c|c}
0&-I\\
\hline
I&0
\end{array}
\right)
=
\left(
\begin{array}{ccc|ccc}
0&i&0&0&&\\
0&0&ib &&-1&\\
0&0&0&&&-1\\
\hline
0&&&0&0&0\\
&0&&-i&0&0\\
&&1&0&-ib&0\end{array}
\right).
\]
Then, with respect to the new, complex symplectic basis 
$$
\tilde X_1:=Y_1-bY_3,\ \tilde X_2:=-b^2X_1+bX_3-iY_2,
\ \tilde Y_1:=-X_1,\ \tilde Y_2:=-iX_2
$$
$$
\tilde X_3:=X_3+bX_1-bX_2+iY_3,\ \tilde Y_3:=\frac i 2(X_3+bX_1+bX_2-iY_3),
$$
$S$ is represented by the block matrix
\[
\tilde S=
\left(
\begin{array}{c|c}
M&0\\
\hline
0&\begin{array}{cc}
-i&0\\
0 & i\end{array}
\end{array}
\right),
\]
where
\[
M:=\left(
\begin{array}{cccc}
0&0&0&0\\
1&0&0&0\\-b^2 & 0&0&-1\\
0 & -(b^2+1)& 0 & 0
\end{array}
\right).
\]
One checks easily that $M$ is 4-step nilpotent, hence conjugate to the matrix
\[
\begin{pmatrix}
0 &&&\\1 & 0 &&\\
&1 & 0 & \\
&&1 & 0
\end{pmatrix}
.
\]
Consequently, $S_r=N_r$ is 4-step nilpotent, and with respect to our complex basis,
$S_r$ is represented by the matrix 
$\tilde S_r=\left(
\begin{array}{c|c}
M&0\\
\hline
0&0
\end{array}
\right)$, and $\tilde S_i$ by 
$\left(
\begin{array}{c|c}
0&0\\
\hline 
0&
	\begin{array}{cc} -i&0\\0&i\end{array}
\end{array}
\right)$. 
SinceÊ\hfill\newline $M\cdot \left(
\begin{array}{cc}
0 & I\\-I & 0\end{array}\right)=
\left(
\begin{array}{cccc}
0 & 0 & 0 & 0\\
0 & 0 & 1 & 0\\
0 & 1 &-b^2&0\\
0 & 0 & 0&-(b^2+1)
\end{array}
\right)$
and
$\begin{pmatrix}
-i&0\\0&i\end{pmatrix}
\begin{pmatrix}
0 & 1\\-1 & 0\end{pmatrix}
=
\begin{pmatrix}0 & -i \\ -i & 0\end{pmatrix}
$, we thus find that 
$$
L_{S_r}= 2 \tilde X_2\tilde Y_1-b^2\tilde Y_1^2-(b^2+1)\tilde
Y_2^2=b^2X^2_1+(b^2+1)X_2^2-2bX_1X_3+2iX_1Y_2
$$
and 
$$
L_{S_i}=-2i\tilde X_3\tilde Y_3
= (X_3+bX_1)^2-b^2X_2^2+Y_3^2+2ibX_2Y_3.
$$
Clearly, neither is $Q_{\Re S_r}\ge 0$, nor $Q_{\Re S_i}\ge 0$.
And, arguing similarly as in the previous example, condition (H) for ($L_{S_r}
+\text{ first order term}$) reduces to solving the system
\[
b^2\xi_1^2+(b^2+1)\xi_2^2-2b\xi_1\xi_3=0, \ \xi_1\eta_2=0,\ \xi_1\xi_2\ne 0.
\]
A solution is given by $\xi_1=\xi_2=1,\ \xi_3=(2b^2+1)/(2b),\ \eta_j=0,\ j=1,2,3$.

In a similar way, one checks that also $L_{S_i}$ satisfies condition (H).

\bigskip

The last two examples show that one can neither dispense with the condition
\eqref{R}, nor  with the condition that $N_r$ be nilpotent of step at most two, in
Proposition \ref{2B}.

We remark that Example \ref{2D} is, in a way, of minimal dimension, if one wants to
show the latter statement. More precisely, it is of minimal possible dimension, if
one requires the nilpotent part $N_r$ to consist of just one Jordan block. This can
easily be  seen from the classification of normal forms of elements in $\sp(n,\CC)$
(see e.g. \cite{hoermander-mehler}, Theorem 2.1), which reveals that nilpotent
elements in
$\sp(n,\CC)$ consisting of just one Jordan block are nilpotent of even step.

\bigskip

\noi {\bf Example \ref{2F}}
On $\HH_2$, consider 
\[
L_S:=(m+c_1)Y_1^2+(m-c_1)Y_2^2+2c_2Y_1Y_2+2i(X_1Y_2-X_2Y_1).
\]
We assume that $c_1,c_2$ 
and $m$ are real, that $c_1^2+c_2^2\ne 0$ and $m\ge\sqrt{c_1^2+c_2^2}$. Then $A$ is
given by 
\[
A=
\left(
\begin{array}{cc|cc}
0&0&0&i\\
0&0&-i&0\\
\hline
0&-i&m+c_1 & c_2\\
i & 0 &c_2 & m-c_1\end{array}
\right)
,
\]
and one verifies readily that $\Re A\ge 0$, since $m\ge\sqrt{c_1^2+c_2^2}$.
Moreover,
$S=D+N$, where $D$ and $N$ are the block matrices
\[
D=\begin{pmatrix}
iJ&0\\C & iJ\end{pmatrix}, \quad N=\begin{pmatrix}
0 & 0\\mI & 0
\end{pmatrix},
\]
with 
$J:=\begin{pmatrix} 0 & 1\\-1 & 0\end{pmatrix}$ and 
$C:=\begin{pmatrix}c_1&c_2\\c_2& -c_1\end{pmatrix}$.

One computes that $D^2=I$, since $CJ+JC=0$.
This implies that $D$ is conjugate to the matrix $\begin{pmatrix}I & 0\\0 & -I\end{pmatrix}$, hence semi-simple. Moreover, $N$ is 2-step nilpotent and commutes with $D$, so that $S=D+N$ is the Jordan decomposition of $S$. Clearly, $\spec S=\{-1,1\}$. But,
\[
\Re Q_D(v)=\trans v\begin{pmatrix}0 & 0\\0 & C\end{pmatrix} v, \quad v\in
\RR^2\times\RR^2,
\]
and $\det C=-(c_1^2+c_2^2)<0$, so that $\Re Q_C$ is an indefinite quadratic form.
This is in agreement with Lemma \ref{3I}, since $D$ is not purely imaginary.

\begin{remark}\label{4E}
{\rm 
In the study  in \cite{mueller-ricci-2} of operators $L_S$ with real matrices $S\in
\sp (n,\RR)$ whose  spectrum is purely real,  it had been most useful to
apply the Jordan decomposition $S=D+N$ of $S$  in order to factorize
$\Gamma_{t,S}^\mu=\Gamma_{t,D}^\mu\times_\mu \Gamma_{t,N}^\mu$. As the previous
example shows, there is no hope of extending this approach to the complex
coefficient case, since $\Re Q_D$ may not be positive semi-definite, so that
$\Gamma_{t,D}^\mu$ may not be tempered.}
\end{remark}


\setcounter{equation}{0}
\section{Twisted convolution and Gaussians \\generated by $\tilde
L^\mu_S$.}\label{gaussians}

Assume that $S\in \sp(n,\CC)$ is such that $\Re Q_S\ge 0$. It is our main goal in
 this section to determine the semigroup generated by the operator
$|U|^{-1}L_S$. Our results  present a generalization of corresponding
results in \cite{mueller-ricci-1},\cite{mueller-ricci-2},
\cite{mueller-ricci-cone}, and are directly related to those in
\cite{hoermander-mehler} by means of the Weyl transform. Instead of
transferring the result from \cite{hoermander-mehler}, Theorem 4.3, by
means of the inverse Weyl transform, we prefer, however,  to give a
direct argument, based on \cite{mueller-ricci-cone}, Theorem 5.2 and ideas
from \cite{hoermander-mehler} and \cite{howe}.

We shall work in the setting of an arbitrary real symplectic vector space
$(V,\sigma)$ of dimension $2n.$ Given two suitable functions $\varphi$ and
$\psi$ on
$V$ and $\mu\in \RR^\times:=\RR\setminus \{0\}$, we define the {\it
$\mu$-twisted convolution} of $\vp$ and $\psi$ as 
\[
\vp \times_\mu \psi (v):=\int_V \vp(v-v')\psi(v')e^{-\pi i\mu
\sigma(v,v')}\,dv',
\]
where $dv'$ stands for the volume form $\sigma^{\wedge (n)}$.

If $f$ is a suitable function on $\HH_V$, we denote by
\[
f^\mu(v):=\int_{-\infty}^\infty f(v,u)e^{-2\pi i\mu u}du
\]
the partial Fourier transform of $f$ in the central variable $u$ at $\mu \in \RR$.

For $\mu\ne 0$, we have 
\begin{equation}\label{5a}
(f*g)^\mu(v)=f^\mu \times_\mu g^\mu (v).
\end{equation}
Moreover, if $A$ is any left-invariant  differential operator on $\HH_V$, then 
there exists a differential operator $\tilde A^\mu$ on $V$ such that 
\begin{equation}\label{5b}
(Af)^\mu =\tilde A^\mu f^\mu.
\end{equation}
Explicitly, if $(x,y)\in \RR^n \times \RR^n$ are coordinates on $V$ associated with a symplectic basis $\{X_j,Y_j\}$, then 
\begin{eqnarray}\label{5c}
\tilde X_j^\mu \vp&=&(\de_{x_j}-\pi i \mu y_j)\vp=\vp \times_\mu(\de
_{x_j}\delta_0),\nonumber\\
\tilde Y_j^\mu \vp&=& (\de_{y_j} +\pi i\mu x_j)\vp=\vp
\times_\mu(\de_{y_j}\delta_0)\\
\tilde U^\mu \vp&=& 2\pi i \mu \vp,\nonumber
\end{eqnarray}
and consequently, 
\begin{equation}\label{5d}
\tilde L_S^\mu=\tilde{\cal L}_A^\mu=\sum_{j,k} a_{jk}\tilde V_j^\mu
\tilde V_k^\mu
\end{equation}
is obtained from $\L_A$ by replacing each $V_j$ in \eqref{1b} by $\tilde
V_j^\mu$. We remark that, for twisted convolutions with $\de_w\delta_0$
on the left, there are analogous formulas:
\begin{eqnarray}\label{5e}
(\de_{x_j}\delta_0)\times_\mu \vp&=& (\de_{x_j}+\pi i \mu y_j)\vp,\\
(\de _{y_j}\delta_0)\times_\mu \vp&=&(\de_{y_j}-\pi i \mu
x_j)\vp.\nonumber
\end{eqnarray}
On $V$, we define the {\it (adapted) Fourier transform} by
\[
\hat f(w):=\int_V f(v)e^{-2\pi i \sigma(w,v)} dv,\quad w\in V.
\]
Observe that then $\hat{\hat f} =f$ and $\int fg=\int \hat f\hat g$, for
suitable functions $f$  and $g$ on $V$.

Consider an arbitrary quadratic form $Q$ on $V^\CC$, with associated Hamilton map $S\in \sp (V^\CC,\sigma)$.
Once we have fixed a symplectic basis $\{X_j,Y_j\}$ of $V$, we may
identify $S$ with a $2n\times 2n$-matrix.
If $\pm \la_1,\dots,\pm \la_m$ are the non-zero eigenvalues of $S$, then 
$\det (\cos S)=\prod_{j=1}^m \cos ^2\la_j$, so that 
the square root 
\begin{equation}\label{5f}
\sqrt{\det (\cos S)}:=\prod_{j=1}^m \cos \la_j
\end{equation}
is well-defined.
Observe that this expression is invariant under all permutations of the
roots of the characteristic polynomial $\det (S-\la I)$, hence an entire
function of the elementary symmetric functions, which are polynomials in 
(the coefficients of) $S$.

Thus, as already observed in \cite{hoermander-mehler},
$\sqrt{\det (\cos S)}$, given by \eqref{5f}, is a well-defined analytic
function of $S\in \sp(V^\CC,\sigma)$.

We shall always consider $T_S:=\tilde L_S^\mu$ as the maximal 
 operator defined by the differential operator
\eqref{5d} on $L^2(V)$; its domain $\dom(T_S)$ consists of all
functions $f\in L^2(V)$ such that $\tilde L_S^\mu f$, defined in the
distributional sense, is in $L^2(V)$.

\begin{lemma}\label{5A}
$\tilde L_S^\mu$ is a closed operator. It is the closure of its restriction to $\S(V)$.
\end{lemma}

The proof of Lemma \ref{5A} will be based on the following well-known 
observation, which follows easily from the formulas \eqref{5c} and
\eqref{5e} (compare also \cite{howe}).

\begin{lemma}\label{5B}
For $w\in V$, denote by $\de_w$ the directional derivative $\de_w f(v)=\frac d{dt}\mid_{t=0} f(v+tw)$, and put 
$\ve_w:=\de_w\delta_0$. Then, the topology in $\S(V)$ is induced by the
semi-norms
\[
||\ve_{w_1} \times_\mu\dots \times_\mu
\ve_{w_N}\times_\mu f\times_\mu\ve_{w_1}\times_\mu\dots \times_\mu
\ve_{w_N}||,
\]
where $w_1,\dots,w_N,w_1',\dots, w_N'$ are arbitrary elements of $V$.
\end{lemma}

\medskip
\noindent{\bf Proof of Lemma \ref{5A}.}
The continuity of $\tilde L_S^\mu$ on ${\cal D}'(V)$ implies the
closedness of the operator $T_S$.

Next, observe that if $f\in L^2(V)$ and $\vp,\psi\in \S(V)$, then it
follows readily from Lemma~\ref{5B} that $\vp \times_\mu f
\times_\mu\psi\in \S(V)$.

Choose a Dirac family $\{\vp_\ve\}$ in $\D(V)$ such that $\vp_\ve(v)=\ve^{-2n} \vp(\ve^{-1}v)$, and assume that $f\in \dom(T_S)$. Then $f_\ve:=\vp_\ve \times_\mu f \times_\mu \vp_\ve \in \S(V)$, and clearly $f_\ve\to f$ in $L^2(V)$ as $\ve\to 0$. Moreover,
\[
\tilde L_S^\mu f_\ve =\vp_\ve \times_\mu f \times_\mu (\tilde L_S^\mu \vp_\ve),
\]
by the left-invariance of $L_S$. And, straight-forward computations based on 
\eqref{5c} -- \eqref{5d} and the symmetry of $L_S$  shows that 
\[
f\times_\mu(\tilde L_S^\mu \vp_\ve)-(\tilde L_S^\mu f)\times_\mu \vp_\ve=f \times_\mu \eta_\ve,
\]
where $\eta_\ve(v)=\ve^{-2n} \eta(\ve^{-1}v),\ \eta\in \D(v)$ and 
$\int \eta \,dv=0$. This implies 
\[
\lim\limits_{\ve\to 0} \tilde L_S^\mu f_\ve=\lim\limits_{\ve\to 0} \vp_\ve \times_\mu(\tilde L_S^\mu f)\times_\mu \vp_\ve=\tilde L_S^\mu f
\]
in $L^2(V)$, and thus $\{f_\ve\}_\ve$ converges in the graph norm to $f.$

\hfill Q.E.D.
\bigskip

Next, obviously the formal adjoint operator of $\tilde L_S^\mu$ is given by
 $\tilde L_{\overline S}^\mu.$ In view of Lemma~\ref{5A}, we thus
have $(\tilde L_S^\mu)^*=\tilde L_{\oS}^\mu$ for the adjoint.

\begin{lemma}\label{5C} Assume that $\Re Q_S\ge 0$. Then the operator
$\tilde \L_A^\mu=\tilde L_S^\mu$ and its adjoint are dissipative, hence
it generates a contraction semigroup $\exp(t\tilde L_S^\mu)$, $t\ge 0$,
on $L^2(V)$.
\end{lemma}

\noindent{\bf Proof.} Clearly, for $f\in \S(V)$, we have 
\begin{eqnarray*}
\Re(\tilde \L_A^\mu f,f)&=&-\Re \sum_{j,k} a_{jk}(\tilde V_j^\mu f,
\tilde V_k^\mu f)\\ &=&-\Re \int_V a_{jk}g_j(v)\overline{g_k(v)} \,dv\\
&=&-\int_V \Re Q_S(g_j(v),\overline{g_k(v)})\, dv\le 0,
\end{eqnarray*}
if we set $g_j:=\tilde V_j^\mu f$. This inequality remains true for arbitrary $f\in \dom (\tilde L_A^\mu)$, by 
Lemma~\ref{5A}, hence $\tilde L_S^\mu$ is dissipative, and the same is true of the adjoint operator $\tilde L_{\overline S}^\mu$, since $\Re Q_S=\Re Q_\oS$.
But then $\tilde L_S^\mu$ generates a contraction semigroup (cf.
\cite{yosida}).

\hfill Q.E.D.
\bigskip

For the case where $\Re Q_S>0$, an explicit formula for the semigroup 
$\exp(\frac{t}{|\mu|} \tilde L_S^\mu)$ has been given in
\cite{mueller-ricci-cone}, Theorem 5.2:

\begin{theorem}\label{5D}
If  $\Re Q_S>0$, then for $f\in L^2(V)$
\begin{equation}\label{5g}
\exp (\frac t{|\mu|} \tilde L_S^\mu)f=f \times_\mu\Gamma_{t,S}^\mu, \quad
t\ge0,
\end{equation}
where, for $t>0$, $\Gamma_{t,S}^\mu$ is a Schwartz function whose Fourier
transform is given by
\begin{equation}\label{5h}
\widehat{\Gamma_{t,S}^\mu}(w)=\frac 1{\sqrt{\det(\cos 2\pi tS)}}
e^{-\frac{2\pi}{|\mu|}\sigma(w,\tan(2\pi tS)w)}.
\end{equation}
\end{theorem}

This result can be extended to the semi-definite case.

\begin{theorem}\label{5E}
Denote by $\sp^+(V^\CC,\sigma)$ the cone of all elements $S\in
\sp(V^\CC,\sigma)$ such that $\Re Q_S\ge 0$. Then the mapping $S\mapsto
\exp(\tilde L_S^\mu)f$ is continuous from $\sp^+(V^\CC,\sigma)$ to 
$L^2(V)$ (respectively to $\S(V)$), if $f\in L^2(V)$ (respectively if
$f\in
\S(V)$), and, for
$S\in \sp^+(V^\CC)$ the mapping $t\mapsto \exp (t \tilde L_S^\mu)f$ is
smooth from $\RR_+$ to $\S(V)$, for every $f\in \S(V)$. Moreover, for
$t\ge 0$, the operator $\exp (\frac t{|\mu|} \tilde L_S^\mu)$ is given by
\eqref{5g},  where $\Gamma_{t,S}^\mu$ is a tempered distribution depending
continuously on $S$, whose Fourier transform is given by \eqref{5h}
whenever $\det(\cos (2\pi t S)\ne 0$.
\end{theorem}

\noindent{\bf Proof.} In order to simplify the notation, let us assume $\mu=1$. We then write $\tilde A$ instead of $\tilde A^1$, 
if $A$ is a left-invariant differential operator on $\HH_V$, and $\vp
\times\psi$ instead of  $\vp\times_1 \psi$.

It is evident from Theorem \ref{5D} that, if $\Re Q_S>0$ and $f\in
\S(V)$, then $f(t)=\exp (t\tilde L_S)f $ is a $C^\infty$-function of $t$
and $S$ with values in $\S(V)$, when $t\ge 0$.

If $w=(w_1,\dots,w_N),w'=(w_1',\dots,w_N')\in V^N,$ we put
$f_{w,w'}:=\ve_{w_1}\times\dots \times \ve_{w_N}\times f \times
\ve_{w_1}\times\dots\times\ve_{w_N'}$.

Then 
\begin{equation}\label{5i}
\frac d{dt} f_{w,w'}=\ve_{w_1}\times\dots \times \ve_{w_N} \times(\tilde
L_Sf)\times\ve_{w_1'}\times\dots \times\ve_{w_N'}.
\end{equation}
But, $\tilde L_S$ commutes with twisted convolutions on the left. Moreover, if $W\in V$ is
 regarded as a left-invariant vector field on $\HH_V$, then for
$j=1,\dots,2n$ we have $ WV_j=V_jW+\sigma(W,V_j)U.$ This implies 
\[
WV_jV_k=V_jV_kW+\sigma(W,V_k)V_jU+\sigma(W,V_j)V_kU,
\]
hence, by some easy computation,
\begin{eqnarray*}
WL_S&=&L_SW +2\sigma(W,\sum\limits_{j,k}a_{jk}V_k)V_jU\\
&=&L_SW+2S(W)U.
\end{eqnarray*}
Taking the partial Fourier transform in the central variable, we obtain 
$$
(\tilde L_S f)\times \ve_W=\tilde W(\tilde L_S f)
=\tilde L_S(\tilde W f)+4\pi i\widetilde{S( W)}f,
$$
i.e.
\begin{equation}\label{5j}
(\tilde L_S f)\times \ve_W=\tilde L_S(f\times \ve_W)+4\pi i f\times 
\ve_{S( W)}.
\end{equation}
Applying this repeatedly to \eqref{5i}, we get 
\begin{equation}\label{5k}
\frac d{dt} f_{w,w'}=\tilde L_{S} f_{w,w'}+4\pi i\sum\limits_{j=1}^N 
f_{w,(w_1',\dots,S(w_j'),\dots,w'_N)}.
\end{equation}
Since $\tilde L_S$ is dissipative, we conclude that 
\begin{eqnarray*}
\frac d{dt} \sum\limits_{w,w'\in \{V_1,\dots,V_{2n}\}^N} ||f_{w,w'}||_2^2
&=&2 \sum\limits_{w,w'\in \{V_1,\dots,V_{2n}\}^N} \Re \left(
\frac{df_{w,w'}}{dt} , f_{w,w'}\right)\\
&\le& C_{N,S}\sum_{w,w'\in \{V_1,\dots,V_{2n}\}^N} ||f_{w,w'}||^2,
\end{eqnarray*}
hence
\begin{equation}\label{5l}
||f(t)||_{(N)}^2:=\sum\limits_{w,w'\in \{V_1,\dots,V_{2n}\}^N} 
||f_{w,w'}||^2\le e^{t C_{N,S}} ||f(0)||_{(N)},
\end{equation}
where $C_{N,S}\le C_N(1+||S||)$.
Notice that, by Lemma \ref{5B}, the semi-norms $||\cdot ||_{(N)}, 
\ N\in \NN$, induce the topology on ${\cal S}(V)$.

Next, assume that $\Re Q_S\ge 0$ and $\Re Q_{S'}>0$, and put $h(t):=\exp(t\tilde L_S)f-\exp (t \tilde L_{S'})f$, where $f\in \S(V)$. Then 
\begin{eqnarray*}
\lefteqn{\frac d{dt}||h(t)||^2}\\
&&=2 \Re (\tilde L_S \exp (t \tilde L_S)f-\tilde L_{S'} \exp (t\tilde L_{S'})f,h(t))\\&&= 2\Re (\tilde L_Sh(t),h(t))+2\Re (\tilde L_{S-S'},\exp (t\tilde L_{S'})f,h(t))\\
&&\le 2C||S-S'||\ ||\exp(t\tilde L_{S'})f||_{(N)}||h(t)||,
\end{eqnarray*}
for some $n\in \NN$ and $C>0$. Thus, because of \eqref{5l},
\[
\frac d{dt}||h(t)||\le C||S-S'||e^{tC_{N,S'}}||f||_{(N)},
\]
hence, if we assume without loss of generality that $C_{N,S'}\ge 1$,
\begin{equation}\label{5m}
||h(t)||\le C||S-S'||(e^{tC_{N,S'}}-1)||f||_{(N)}.
\end{equation}
From \eqref{5m} one deduces that $\exp(t\tilde L_S)f$ is continuous as a
function of $S$ with values in $L^2(V)$, first, for $f\in {\cal S}(V)$,
but then also for arbitrary $f\in L^2(V)$, by the contraction property.
Once this is shown, it follows easily with the aid of \eqref{5l} that
$\exp(t\tilde L_S)f$ is also continuous  as a function of $S$ with values
in $\S(V)$, given $f\in \S(V)$. In particular, if $S=\lim S'$, with $\Re
Q_{S'}>0$, then $\Gamma_{t,S'}^1$, given by \eqref{5h}, converges in
$\S'(V)$ towards a tempered distribution $\Gamma_{t,S}^1$, so that
\eqref{5g} holds true for arbitrary $S\in \sp^+(V^\CC,\sigma)$. 

Moreover,
the dominated convergence theorem shows that also formula \eqref{5h}
remains valid whenever $\det(\cos 2\pi tS)\ne 0$. Clearly, also the
mapping $S\mapsto \Gamma_{t,S}^1\in \S'(V)$ is continuous.

 Finally, since
\[
\S(V)\subset \dom(\tilde L_S^k) \quad\mbox{ for every } k\in \NN,
\]
it follows easily from \eqref{5i} and \eqref{5l} that the mapping
$t\mapsto
\exp(t\tilde L_S)f$ is smooth from $\RR_+$ to $\S(V)$ if $f\in \S(V)$.

\hfill Q.E.D.
\bigskip

Observe that Theorem \ref{5E} implies that 
\begin{equation}\label{5n}
\Re \sigma (w,\tan (2\pi tS)w)\ge 0 \quad \forall w\in V, \ t\ge 0,
\end{equation}
whenever $\det (\cos 2\pi tS)\ne 0$.

In the coordinates $v=(x,y)\in \RR^n\times \RR^n$, the symplectic Fourier
 transform can be written as
\[
\hat f(w)=\int f(v)e^{-2\pi i \trans w Jv} \,dv,
\]
and one computes that $(\de_{ v_j} f)^\wedge(w) =-2\pi i (J w)_j\hat f(w),
\ ((J v)_j f)^\wedge(w) =-\frac 1 {2\pi i} \de_{w_j}\hat f(w)$. This shows
that
$(\tilde V_j ^\mu f)^\wedge=\hat V_j^\mu \hat f$, where explicitly 
\begin{equation}\label{5o}
\hat V_j^\mu =\frac \mu 2 \de_{ w_j}-(2\pi i)(J w)_j,\quad j=1,\dots,2n.
\end{equation}
Of course, $(\tilde L_S^\mu f)^\wedge =:\hat L_S^\mu \hat f$, where 
\begin{equation}\label{5p}
\hat L_S^\mu =\sum_{j,k} a_{jk} \hat V_j^\mu \hat V_k^\mu.
\end{equation}
If $\Re Q_S\ge 0$, then it follows from Theorem 5.5 that $|\mu| \frac \de{\de t} \Gamma_{t,S}^\mu =\tilde L_S^\mu \Gamma_{t,S}^\mu$ in the sense 
of distributions. Taking Fourier transforms, it is clear by \eqref{5h}
that the corresponding formula for the Fourier transforms will also hold
pointwise, i.e. 
\begin{equation}\label{5q}
|\mu|\frac \de{\de t} \widehat{\Gamma_{t,S}^\mu}(w)=\tilde L_S^\mu \widehat{\Gamma_{t,S}^\mu}(w)\quad \forall w\in V,
\end{equation}
whenever $\det (\cos 2\pi tS)\ne 0$.

For {\it arbitrary} $S\in \sp (V^\CC,\sigma)$, and 
complex $t\in \CC,\ w\in V^\CC$, let us {\it define}
$\widehat{\Gamma_{t,S}^\mu}(w)$ by formula \eqref{5h}, whenever $\det(\cos
2\pi tS)\ne 0$. Observe that $\widehat{\Gamma_{t,S}^\mu }$ may not be
tempered, if 
$
S\not\in  \sp^+ (V^\CC,\sigma)\ \text{or}\ t\not\in \RR_+.
$
By analytic extension, formula \eqref{5q} then remains valid, i.e.
\begin{equation}\label{5r}
|\mu|\frac\de{\de t} \widehat{\Gamma_{t,S}^\mu} (w)=\tilde L_S^\mu
\widehat{\Gamma_{t,S}^\mu} (w)\quad \forall w\in V^\CC, t\in \CC, S\in
\sp(V^\CC,\sigma),
\end{equation}
whenever $\det(\cos 2\pi tS)\ne 0$, if we denote by $\frac\de{\de t}$ 
and $\de_{w_j}$ the complex derivatives with respect  to $t\in
\CC$ and the complex variable $w_j$ in \eqref{5o}, respectively.

This allows us to introduce complex symplectic changes of coordinates. Let 
$T=(T_{jk}) \in \Sp(n,\CC)$ be an arbitrary symplectic matrix, and
introduce new symplectic coordinates
\[
z=Tw\in \CC^{2n}, \ w\in \RR^{2n}.
\]
Since $\trans TJT=J$, we have 
\[
Jw=\trans T(Jz),
\]
hence
\[
\hat V_j^\mu =\frac \mu 2 \sum_k T_{kj}\de _{z_k}-2\pi i \sum_k
T_{kj}(Jz)_k,
\]
when acting on holomorphic functions
 (such as $\widehat{\Gamma_{t,S}^\mu}$). Putting, in analogy with
\eqref{5o},
\[
Z_k:=\frac\mu 2 \de_{z_k}-2\pi i (Jz)_k,\quad k=1,\dots 2n,
\]
 we get
\begin{equation}\label{5s}
\hat V_j^\mu =\sum_k T_{kj} Z_k.
\end{equation}
Assume now that  we have a splitting of the (complex) symplectic
coordinates in two blocks, $z'=(z_1,\dots,z_q;z_{n+1},\dots z_{n+q})$ and
$z''=(z_{q+1},\dots, z_{n}; z_{n+q},\dots,z_{2n})$, where $1\le q <n$.
Then, the following lemma is obvious.

\begin{lemma}\label{5F} 
Let $f(z)=f_1(z')f_2(z'')$, where $f_1$ and $f_2$ are holomorphic functions. Then, for $k\in \{1,\dots,q\}\cup \{n+1,\dots,n+q\},\ Z_kf_1$ is again a function of $z'$, and 
\[
Z_k f(z)=(Z_kf_1)(z')f_2(z'').
\]
\end{lemma}
For $S\in \sp(V^\CC,\sigma)$, denote again by $S=S_r+S_i$ the
decomposition given by (2.12).

\begin{proposition}\label{5G}
Assume that $\det \cos (2\pi tS)\ne 0$. 
Then
\begin{equation}\label{5t}
\widehat{\Gamma_{t,S}^\mu} (w)=\widehat{\Gamma_{t,S_r}^\mu}(w)\widehat{\Gamma_{t,S_i}^\mu} (w),
\end{equation}
and
\begin{eqnarray}\label{5u}
\hat L_{S_r}^\mu \widehat{\Gamma_{t,S}^\mu}(w)=(\hat L_{S_r}^\mu
\widehat{\Gamma_{t,S_r}^\mu)}(w)\widehat{\Gamma_{t,S_i}^\mu}(w)
=|\mu|(\de_t\widehat{\Gamma_{t,S_r}^\mu})(w)\widehat{\Gamma_{t,S_i}^\mu}(w),
\end{eqnarray}
if $t\in \CC, \ w\in V^\CC$.
\end{proposition}

\bigskip

\noindent{\bf Proof.} Choosing real symplectic coordinates, we may assume that $w\in \RR^{2n}$. Since $S_r$ and $S_i$ correspond
 to different sets of Jordan blocks of $S$, we can choose $T\in
\Sp(n,\CC)$ such that
\[
\tilde S:=TST^{-1}=\left(
	\begin{array}{cc}
		\tilde S_r & 0\\
		0             &\tilde S_i
	\end{array}
	\right)
\]
with respect to suitable blocks $z'$ and $z''$ of complex symplectic
coordinates, say
$z'\in \RR^{2q},\ z''\in \RR^{2(n-q)}$, where 
\[
TS_rT^{-1}=
\left(
	\begin{array}{cc}
		\tilde S_r & 0\\
		0             &0
	\end{array}
	\right),\quad 
TS_iT^{-1}=\left(
	\begin{array}{cc}
		0& 0\\
		0             &\tilde S_i
	\end{array}
	\right).
\]
In the new coordinates $z=Tw,\ \widehat{\Gamma_{t,S}^\mu}$ has the form
\begin{eqnarray}\label{5v}
\widehat{\Gamma_{t,S}^\mu}&=&\frac 1
	{\sqrt{\det(\cos 2\pi tS)}} e^{-\frac{2\pi}{|\mu|}\trans z\cdot J\tan (2\pi t\tilde S)\cdot z}\\\nonumber
&=&\widehat{\Gamma_{t,\tilde S_r}^\mu}(z')\widehat{\Gamma_{t,\tilde S_i}^\mu}(z''),
\end{eqnarray}
which proves \eqref{5t}. Moreover,
\[
\hat L_{S_r}^\mu =\sum_{j,k} b_{jk} \hat V_j^\mu \hat V_k^\mu,
\]
where $B=(b_{jk})=S_rJ$, hence by \eqref{5s},
\begin{eqnarray*}
\hat L_{S_r}^\mu&=& \sum_{j,k}\sum_{l,m} b_{jk} T_{lj}Z_lT_{mk}Z_m\\
&=& \sum_{l,m}(TB\trans
T)_{lm}Z_lZ_m=\sum_{l,m}\Big ((TS_rT^{-1})J\Big )_{lm}Z_lZ_m.
\end{eqnarray*}
Putting $C:=\tilde S_r J$, which is an $2q\times 2q$-matrix, we find that 
\[
\hat L_{S_r}^\mu =\sum_{j,k=1}^{2q} C_{jk}Z_jZ_k.
\]
Formulas \eqref{5u} are now a consequence of \eqref{5v}, Lemma~\ref{5F}
and \eqref{5r}

\hfill Q.E.D.


\setcounter{equation}{0}
\section{Reduction to Hamiltonians with purely \\real spectrum}
\label{reduction}

In this section, we shall prove Theorem 2.1. So, assume that $\Re Q_S\ge
0$ and $S_i\ne 0$, and that (R) holds. We begin with the
case where $|\Re \al|<\nu$. The following proposition will imply Theorem
2.1 (i).

\begin{proposition}\label{6A}
For $f\in \S(\HH_V)$, the integral 
\begin{equation}\label{6a}
\l K_\al, f\r:=-\int_{-\infty}^{-\infty} \ \int_{0}^{+\infty}
e^{-2\pi\al t\,\sgn\mu} \langle \Gamma_{t,S}^\mu , f^{-\mu} \r\,
dt\,\frac{d\mu}{|\mu|}
\end{equation}
converges absolutely and defines a tempered distribution $K_\al$ for $|\Re \al| < \nu$. 
Moreover, $K_\al$ is a fundamental solution for $L_{S,\al}$, i.e.
$L_{S,\al} K_\al=\delta _0$.
\end{proposition}
Here, $\Gamma_{t,S}^\mu\in \S'(V)$ is given by Theorem \ref{5E}.

\bigskip

\noindent{\bf Proof.}
Recall that $S=S_r+S_i$, where 
\[
\spec S_i=\{\pm \om_1,\dots,\pm \om_{n_1}\} \subset \CC\setminus \RR
\]
and $\nu_j=\Im \om_j>0,\  j=1,\dots,n_1$. Also, $\nu
=\sum\limits_{j=1}^{n_1}
\nu_j$, $\nu_{min}=\min\limits_ {j=1,\dots,n_1} \nu_j$. Let
\[
\spec S_r\setminus \{0\}=\{\pm \la_1,\dots,\pm \la_{n_2}\},
\]
where $\la_k>0, \ k=1,\dots,n_2$. Then, by Theorem \ref{5E},
\begin{equation}\label{6b}
\widehat{\Gamma_
	{\frac t{2\pi},S}^\mu}(w)=\frac 1
	{\sqrt{\det(\cos tS)}} e^{-\frac{2\pi}{|\mu|}\sigma(w,\tan (tS)w)},
\end{equation}
whenever
\begin{equation}\label{6c}
\sqrt{\det(\cos tS)}=\prod_{k=1}^{n_2}\cos (t\la_k)
\prod_{j=1}^{n_1}\cos(t\om_j)\ne 0.
\end{equation}
Thus, potential singularities in \eqref{6b} arize when $t\la_k=\frac \pi
2+\ell \pi$ for some $k\in \{1,\dots,n_2\}$ and $\ell\in \ZZ$. By means
of partial Fourier transforms, we shall show that these points are in
fact not singular, if we consider $t\mapsto \widehat{\Gamma_{\frac
t{2\pi},S}^\mu}$ as a family of distributions.

For any subset $I\subset \{1,\dots,n_r\}$, denote by $V_I \subset V$ the real subspace
\[
V_I:=\sum_{k\in I}^{\quad\oplus} V_{\la_k}\oplus V_{-\la_k}.
\]
Then $V_I$ and $V_I^\bot$ are $S_1$ and $S_2$-invariant symplectic subspaces.
If we choose real symplectic coordinates $w'$ for $V_I$ and $w''$ for $V_I^\bot$, then $S$ will be represented by a block matrix
\[
S=\left(
\begin{array}{cc}
S_I & 0\\
0 & S_{I^\bot}\end{array}\right)
\]
with respect to the coordinates $(w',w'')$ for $V$,
and similarly as in Section \ref{gaussians} we find that 
\begin{equation}\label{6d}
\widehat{\Gamma_{t,S}^\mu} (w',w'')=\widehat{\Gamma_{t,S_I}^\mu}(w')
\widehat{\Gamma_{t,S_{I^\bot}}^\mu}(w''),
\end{equation}
where
\begin{equation}\label{6e}
\widehat{\Gamma_{\frac t{2\pi},S_I}^\mu}(w')=
\frac{1}
{\prod_{k\in
I}\cos(t\la_k)}e^{-\frac{2\pi}{|\mu|}\sigma(w',\tan(tS_I)w')},
\end{equation}
\begin{equation}\label{6f}
\widehat{\Gamma_{\frac t{2\pi},S_{I^\bot}}^\mu}
(w'')=\frac 1{\prod_{k\not\in I} \cos (t\la_k)\prod_j\cos
(t\om_j)}e^{-\frac{2\pi}{|\mu|}\sigma(w'',\tan(tS_{I^\bot})w'')}.
\end{equation}
From \eqref{6e} one computes that
\begin{equation}\label{6g}
\Gamma_{\frac t{2\pi},S_I}^\mu (v')=
	\frac{c|\mu|^{|I|}}{\prod_{k\in I} \sin(t\la_k)} 
	e^{-\frac \pi 2|\mu|\sigma(v',\cot(tS_I)v')},
\end{equation}
where $c$ is a constant of modulus 1 (see e.g. \cite{hoermander-mehler}, 
Theorem 7.6.1). This
 identity holds, unless $t\la_k=\ell \pi$ for some $k\in I$ and $\ell\in
\ZZ$. Observe also that all exponentials in \eqref{6e} -- \eqref{6g} are
bounded by 1.

\begin{lemma}\label{6B}
There exist a constant $C>0$ and a Schwartz norm $||\cdot ||_\S$ on
$\S(V)$, such that, for every 
$0\le \delta\le 1$,
\begin{equation}\label{6h}
\left |\l
	\Gamma_{\frac t{2\pi},S}^\mu,\vp
\r\right |
\mid 
\le C 
\frac
	{|\mu|^{\delta n_2}(1+|\mu|^{n_2})^{1-\delta}}
{\prod\limits_{k=1}^{n_2}|\sin(t\la_j)|^\delta
\prod\limits_{j=1}^{n_1}|(\cos (t\om_j)|} \,||\vp||_\S, \quad \forall t\ge
0,
\end{equation}
for every $\vp\in \S(V)$.
\end{lemma}

\noindent{\bf Proof.} Given $t\ge 0$, define a subset $I=I_t$ of
$\{1,\dots,n_2\}$ as follows:

\noindent $k\in \{1,\dots,n_2\}$ belongs to $I$ if
and only  if there exists an $\ell \in \ZZ$ such that $|t\la_k-\frac \pi
2 -\ell\pi|\le \frac \pi 4$. Then, for $k\in I$ we have
$|\sin(t\la_k)|\ge\cos \frac \pi 4>0$, and if $k\not\in I$, then
$|\cos(t\la_k)|\ge \cos \frac \pi 4$. Since 
\begin{eqnarray}\label{6i}
\l \Gamma_{t,S}^\mu,\vp\r&=&\l\widehat{\Gamma_{t,S}^\mu},\hat
\vp\r=\l\widehat{\Gamma_{t,S_I}^\mu}\otimes
\widehat{\Gamma_{t,S_{I^\bot}}^\mu},\hat\vp\r\\
	&=&\l\Gamma_{t,S_I}^\mu\otimes \widehat{\Gamma_{t,S_{I^\bot}}^\mu},\hat
\vp^I\r,\nonumber
\end{eqnarray}
where $\hat\vp^I$ denotes the partial Fouriertransform in $v'$, the
formulas \eqref{6f} and \eqref{6g} therefore imply 
\[
\left|\l\Gamma_{t,S}^\mu,\vp\r\right |\le \frac
	{C|\mu|^{|I|}}
	{\prod_j|\cos(t\om_j)|}\,||\hat\vp^I||_1.
\]
This gives \eqref{6h} for $\delta=0$.

On the other hand, choosing $I=\{1,\dots,n_2\}$, we have
\[
\left |\l\Gamma_{\frac t{2\pi},S}^\mu,\vp\r\right |\le C\frac
	{|\mu|^{n_2}}
	{\prod_{k=1}^{n_2} |\sin(t\la_k)|\prod_{j=1}^{n_1}|\cos(t\om_j)|} \,
	||\vp||_\S,
\]
which is the case $\delta=1$, and \eqref{6h} follows immediately from
these extreme cases by interpolation.

\hfill Q.E.D.
\bigskip

\noindent Observe next that for $t>0$
\begin{equation}\label{6j}
\frac 1{\cos (t\om_j)} = \frac{2e^{it\om_j}}
	{1+e^{2it\om_j}}=2\sum_{m=0}^\infty (-1)^me^{(2m+1)it\om_j},
\end{equation}
which implies 
\begin{equation}\label{6k}
\frac 1
	{\prod\limits_j|\cos t\om_j|} =O(e^{-t\nu}),\quad t\ge 0.
\end{equation}
The integral in \eqref{6a} can thus be estimated in modulus by
\[
C\int\limits_{-\infty}^\infty \ \int\limits_0^\infty \frac
	{e^{-2\pi(\nu-|\Re\al|)t}}
	{\prod\limits_{j=1}^{n_2} |\sin(t\la_j)|^\delta}\, dt
	|\mu|^{\delta n_2-1}(1+|\mu|^{n_2})^{1-\delta} ||f^{-\mu}||_\S \, d\mu,
\]
which is convergent if we choose $\delta>0$ sufficiently small, provided
that $|\Re \al|<\nu$. One also checks easily that
$L_{S,\al}K_\al=\delta_0$ (compare the proof of
Theorem 6.1 in  \cite{mueller-ricci-cone}). This completes the proof of
Proposition \ref{6A}.

\hfill Q.E.D.
\medskip

\begin{remark} \label{rem6}
{\rm If we argue in a similar way  in Example \ref{2C}, by choosing, for a
given $t>0,$ either the expression for $\widehat{\Gamma_{t,S}^\mu}$ or for
$\Gamma_{t,S}^\mu$ in order to carry out the estimations, we find that 
the statement of Proposition \ref{6A} remains true for Example \ref{2C}.
This shows that Proposition \ref{6A} may be true even when
property \eqref{R} is not satisfied.}

\end{remark}
\medskip

\noindent In order to prove Theorem 2.1(ii), let us put
\[
R:=L_{S_r},\ R_\beta:=L_{S_r}+i\beta U,\quad \beta\in \CC.
\]
Assume that $\beta_1,\dots,\beta_N$ are analytic functions of $\al$, and
set, for $|\Re\al|<\nu$,
\begin{equation}\label{l}
K_\al^N:=UR_{\beta_N(\al)}\dots R_{\beta_1(\al)}K_\al\in \S'(\HH_V).
\end{equation}

\begin{lemma}\label{6C}
Assume that $\beta_1,\dots,\beta_N$ have been chosen in such a way that
 the family of tempered distributions $K_\al^N$ extends analytically
from the strip $|\Re\al|<\nu$ to the wider strip $|\Re \al|<M$. Then
$L_{S,\al}$ is locally solvable for $|\Re \al|<M$, provided the operators
$R_{\beta_j(\al)}, j=1,\dots,N$, are locally solvable.
\end{lemma}

\noindent{\bf Proof.} As $[S,S_r]=0$, all operators $L_{S,\al}, R_{\pm
\beta_j}$ and $U$ commute. Thus, if $|\Re \al|< \nu$,
\begin{eqnarray}\label{6m}
L_{S,\al}K_\al^N&=&UR_{\beta_N}\dots R_{\beta_1}L_{S,\al} K_\al\\
&=&UR_{\beta_N}\dots R_{\beta_1}\delta_0,\nonumber
\end{eqnarray}
where the $\beta_j$ have to be evaluated at $\al$.
By analyticity, this identity remains valid for $|\Re\al|<M$.

Since $U,R_{\beta_1(\al)},\dots,R_{\beta_N(\al)} $ are locally solvable,
 this implies local solvability of $L_{S,\al}$ (compare the proof of
 Lemma 7.4 in \cite{mueller-ricci-cone}).

\hfill Q.E.D.
\bigskip

\noindent Let us examine when the family of distributions $K_\al^N$, 
$|\Re \al|<\nu$, can be extended analytically to a wider strip. By
\eqref{6a}
\begin{equation}\label{6n}
\l K_\al^N,f\r
	=-2\pi i\int_{-\infty}^{+\infty} \ \int_0^{+\infty} e^{-2\pi\al
t\,\sgn\mu}
	\l\tilde R_{\beta_N}^\mu\dots \tilde R_{\beta_1}^\mu
\Gamma_{t,S}^\mu,f^{-\mu}\r\, dt \, \sgn\mu\, d\mu.
\end{equation}
We decompose the integration in $\mu$ into $\int_0^{+\infty}d\mu$ and $\int_{-\infty}^0 d\mu$, and denote 
the corresponding contributions by $\l K_\al^+,f\r$ and $\l
K_\al^-,f\r$, respectively. In the sequel, we shall only regard
$K_\al^+$, since the discussion of $K_\al^-$ is similar. So, assume that
$\mu>0$. Observe that 
\begin{equation}\label{6o}
\l K_\al^+,f\r
	=-2\pi i\int_0^{+\infty}\ \int_0^{+\infty} e^{-2\pi \al t}
	\l \hat R_{\beta_N}^\mu\dots \hat R_{\beta_1}^\mu \
\widehat{\Gamma_{t,S}^\mu},
\widehat{f^{-\mu}}\r\, dt\, d\mu.
\end{equation}
Moreover, by Proposition \ref{5G},
\begin{eqnarray*}
\lefteqn{e^{-2\pi \alpha t}
	\l \hat R_{\beta_N}^\mu\dots \hat R_{\beta_1}^\mu
\widehat{\Gamma_{t,S}^\mu},\, \vp\r}\\
	&&=e^{-2\pi \al t}
	\l\hat R_{\beta_N}^\mu\dots \hat R_{\beta_2}^\mu [(\hat
L_{S_r}^\mu-2\pi\beta_1\mu)\widehat{\Gamma_{t,S_r}^\mu} ]\widehat{
\Gamma_{t,S_i}^\mu},\, \vp\r\\ &&=\mu\l\hat R_{\beta_N}^\mu\dots \hat
R_{\beta_2}^\mu[\de _t(e^{-2\pi\beta_1t}\widehat
{\Gamma_{t,S_r}^\mu})][e^{2\pi(\beta_1-\al)t}\widehat{\Gamma_{t,S_i}^\mu}],
\, \vp\r,
\end{eqnarray*}
hence
\begin{eqnarray}\label{6p}
\lefteqn
	{e^{-2\pi\al t}\l\hat R_{\beta_N}^\mu\dots \hat
R_{\beta_1}^\mu\widehat{\Gamma_{t,S}^\mu},\,
	\widehat{f^{-\mu}}\r} \\\nonumber
&&=\mu\frac d{dt}
	\l\hat R_{\beta_N}^\mu\dots \hat R_{\beta_2}^\mu e^{-2\pi\al t}
	\widehat{\Gamma_{t,S}^\mu},\, \widehat{f^{-\mu}}\r\\
&&\qquad -\mu e^{-2\pi \beta_1 t} \l\hat R_{\beta_N}^\mu \dots 
\hat R_{\beta_2}^\mu(\widehat{\Gamma_{t,S_r}^\mu }\
\de_t[e^{2\pi(\beta_1-\al)t}
 \widehat{\Gamma_{t,S_i}^\mu}]),\,\widehat{f^{-\mu}}\r.\nonumber
\end{eqnarray}
Inserting this into \eqref{6o} and integrating by parts (observe that all
expressions are smooth in $t\ge 0$, by Theorem \ref{5E}), we obtain
\begin{eqnarray}\label{6q}
\nonumber
\l K_\al^+,f\r&=& 2\pi i\int_0^{+\infty} \l\hat R_{\beta_N}^\mu\dots
\hat R_{\beta_2}^\mu \1,\widehat{f^{-\mu}}\r \, \mu\, d\mu\\ 
	&+&2\pi i\int_0^{+\infty} \int_0^{+\infty} e^{-2\pi\beta_1 t}\left\l \hat
R_{\beta_N}^\mu \dots \hat R_{\beta_2}^\mu(\widehat
{\Gamma_{t,S_r}^\mu}\,
\de_t[e^{2\pi(\beta_1-\al)t}\widehat{\Gamma_{t,S_i}^\mu}]),\widehat{f^{-\mu}}
\right \r
dt\mu \,d\mu.
\end{eqnarray}
Notice that the boundary term at $t=+\infty$ vanishes, because of
\eqref{6h} and \eqref{6k}, since $|\Re\al |<\nu$.

Now, fix some number $M>\nu$, and denote by 
\[
\vp(t):=\frac
	1
{\sqrt{\det \cos (tS_i)}}=
\frac
	1
{\prod\limits_{j=1}^{n_1} \cos (t\om_j)}.
\]
By \eqref{6j}, $\vp(t)$ has the series expansion 
\[
\vp(t)=2^{n_1}\sum_{m\in \NN^{n_1}} (-1)^{|m|}e^{it\sum_j(2m_j+1)\om_j},
\quad t>0.
\]
We truncate this series to those $m$ for which 
$\sum_j(2m_j+1)\nu_j<M$, so that, for $t>0,$
\begin{eqnarray}\label{6r}
\vp(t)=2^{n_1}\sum_{\sum_j(2m_j+1)\nu_j<M} e^{it\sum_j(2m_j+1)\om_j}+
	O(e^{-Mt}).
\end{eqnarray}
Observe that $\widehat{\Gamma_{t,S}^\mu}:=\widehat{\Gamma_{t,S_r}^\mu}\ \widehat{\Gamma_{t,S_i}^\mu}$, where
\begin{equation}\label{6s}
\widehat{\Gamma_{\frac t{2\pi},S_i}^\mu}(w)=\vp(t)e^{-\frac{2\pi}{|\mu|}\sigma(w,\tan(tS_i)w)}.
\end{equation}
According to \eqref{6r}, we split $\vp(t)$ into a finite sum of terms, and
correspondingly the integral \eqref{6o} into a finite sum of integrals.

Arguing as in the proof of Proposition 6.1, the integral containing the remainder term is absolutely convergent  for $|\Re\al|<M$ and depends analytically on $\al $ in this region.

We then only have to discuss the other terms, and to this end, we imagine 
that $\vp(t)$ in \eqref{6s} has been replaced by
\begin{equation}\label{6t}
\tilde \vp(t)=e^{it\sum_j(2m_j+1)\om_j},
\end{equation}
for some $m\in \NN^{n_1}$ such that $\sum_j(2m_j+1)\nu_j<M$, i.e.
that $\widehat{\Gamma_{\frac t {2\pi},S_i}^\mu}$ has been replaced in
 \eqref{6p} by $\tilde \vp(t)e^{-\frac{2\pi}{|\mu|}Q_t}$, where
\[
Q_t(w):=\sigma(w,\tan (tS_i)w).
\]
We choose, with $m$ as in \eqref{6t},
\begin{equation}\label{6u}
\beta_1:=\al-\sum_i(2m_j+1)i\om_j.
\end{equation}
Then $\Re\beta_1=\Re \al-\sum_j(2m_j+1)\nu_j$, and 
\[
e^{(\beta_1-\al)t}\tilde \vp(t)e^{-\frac{2\pi}\mu Q_t}=e^{-\frac{2\pi}\mu Q_t}.
\]
The next lemma can be proved along the same lines as Lemma 6.3 in
\cite{mueller-ricci-cone}.

\begin{lemma}\label{6D}
There exist  quadratic forms $Q_{jk}$ on $V$ such that
\[
Q_t(w)=\sigma(w,\tan(tS_i)w)=\sum_{j=1}^{n_1}\sum_{k=0}^\ell
t^k\tan^{(k)}(\om_j t)Q_{jk}(w),
\]
where $\tan^{(k)}$ is the $k$-th derivative of the tangent function and $\ell+1$ is the dimension of the largest Jordan block of $S_i$.
\end{lemma}
We now obtain 
\[
\frac{\de Q_t}{\de t}(w)=\sum_{j=1}^{n_1} \sum_{k=0}^\ell t^k\tan ^{(k+1)}
(\om_j t)\tilde Q_{jk}(w),
\]
for some other quadratic forms $\tilde Q_{jk}$. Since 
\[
\de _t(e^{-\frac {2\pi} \mu Q_t(w)})=
	-\frac{2\pi}\mu \ \frac{\de Q_t}{\de t}(w)e^{-\frac{2\pi}\mu Q_t(w)},
\]
the second term in \eqref{6p} then decomposes as a finite sum of terms of
the form
\begin{eqnarray}\label{6v}
\nonumber
c\int_0^{+\infty}\ \int_0^{+\infty} &&e^{-2\pi\beta_1t} 
\tan ^{(k+1)} (2\pi \om_{j_0}t)t^k \\
 &&\cdot\left\l 
\hat R_{\beta_N}^\mu\dots \hat R_{\beta_2}^\mu
\left(
	\tilde Q_{j_0k}\widehat{\Gamma_{t,S_r}^\mu}e^{-\frac{2\pi}\mu Q_t}\right) ,
\widehat{f^{-\mu}}\right\r \,dt\, d\mu.
\end{eqnarray}
It is now important to make the following observation:

\medskip
If we choose complex symplectic coordinates $z'$ and $z''$ corresponding
to Jordan blocks of $S_r$ and $S_i$, respectively, as in the proof of
Proposition \ref{5G}, then $\widehat{\Gamma_{t,S_r}^\mu}$ depends only on
$z'$, $Q_t$ and $\tilde Q_{jk}$ on $z''$, and $\hat
R_{\beta_j}^\mu=L_{S_r}-2\pi\beta_j\mu$ on $z'$. Therefore, 
\begin{eqnarray*}
\lefteqn
{\left\l \hat R_{\beta_N}^\mu \dots \hat R_{\beta_2}^\mu(\tilde
Q_{j_0k}\widehat{\Gamma_{t,S_r}^\mu}e^ {-\frac{2\pi}\mu Q_t}),\,\widehat
{f^{-\mu}}\right\r}\\ &&
=\frac 1 {\vp(2\pi t)} \left\l
	\tilde Q_{j_0 l} \hat R_{\beta_N}^\mu\dots \hat R_{\beta_3}^\mu
([(\hat
L_{S_r}^\mu-2\pi\beta_2\mu)\widehat{\Gamma_{t,S_r}^\mu}]
\widehat{\Gamma_{t,S_i}^\mu}),
\widehat{f^{-\mu}}\right\r.
\end{eqnarray*}
Moreover,
\begin{equation}\label{6w}
\tan(\om_jt)=i-2i\sum_{p=0}^\infty (-1)^pe^{2i(p+1)\om_jt}, \quad t>0,
\end{equation}
so that 
\[
\tan^{(k+1)}(\om_jt)=-2i\sum_{p=0}^\infty
(-1)^p(2i(p+1))^{k+1}e^{2i(p+1)\om_jt}.
\]
This shows that 
\begin{equation}\label{6x}
\cos (\om_{j_0}t)\tan ^{(k+1)}(\om_{j_0}t)=\sum_{p=0}^\infty
a_pe^{it(2p+1)\om_{j_0}t},\quad t>0
\end{equation}
and then 
\[
e^{-2\pi\beta_1t}\, \frac 1
{\vp(2\pi t)}\tan^{(k+1)}(2\pi \om_{j_0}t)=e^{-2\pi(\al-2i\om_{j_0}t)}
\sum_{q\in \NN^n} b_qe^{2\pi it\,2q_j\om_jt},
\]
where the latter series converges locally absolutely and uniformly for
$t>0$.

If we truncate this series in a similar way as before, such that the remainder term is analytic in 
$|\Re \al|<M$, and put everything together, we see that each term
\eqref{6v} can be further decomposed into terms which, except for the
remainder term, are of the form
\begin{eqnarray}\label{6y}
\nonumber
\lefteqn{c\int_0^{+\infty} \ \int_0^{+\infty}e^{-2\pi 
(\al-\sum_j 2q_ji\om_j)t}\, t^k}\\
&&\cdot\left\l \tilde Q_{j_0k} \hat R_{\beta_N}^\mu\dots \hat
R_{\beta_3}^\mu([\hat R_{\beta_2}^\mu
\widehat{\Gamma_{t,S_r}^\mu}]\widehat{\Gamma_{t,S_i}^\mu},
\widehat{f^{-\mu}}\right\r\, dt\,d\mu,
\end{eqnarray}
where $q=(q_1,\dots,q_{n_1})\in \NN^{n_1}$ and $q_{j_0}\ge 1$, i.e.
$|q|\ge 1$.

By Lemma \ref{6B}, this integral converges absolutely for
 $\Re \al > -(\nu+2\nu_{j_0}+\sum_j 2q_j\nu_j)$ and not only for $\Re \al
>-\nu$. Thus, the contribution to $K_\al^N$ given by the
integration over $\mu>0$ has been extended as an analytic family of
distributions from
$\Re \al >-\nu$ to $\Re \al >-\nu-2\nu_{\min}$.

We can at this point iterate the argument above, in order to extend
$K_\al^+$ to the domain $\Re \al >-M$. If one compares \eqref{6y} with
\eqref{6o}, one finds that the only new features are the presence of the
quadratic forms $\tilde Q_{j,k}$ and the powers $t^k$ of $t.$ The factors 
$\tilde Q_{j,k}$ are harmless, since they only depend on $z'',$ so that
the multiplication with $\tilde Q_{j,k}$ commutes with each of the
operators $\hat R_{\beta_j}^\mu.$

As for powers $t^k$ of $t,$ observe that 
\begin{eqnarray}\label{6z}
&&t^k (\hat R_{\beta_2}^\mu)^{k+1}\widehat{\Gamma_{t,S_r}^\mu}\\
&&\quad=\mu\,e^{2\pi\beta_2 t}\frac \de{\de t}
\left [\sum_{j=0}^k (-1)^j\frac{k!}{(k-j)!}(t\hat R_{\beta_2}^\mu)^{k-j}
(e^{-2\pi\beta_j t}\widehat{\Gamma_{t,S_r}^\mu})\right ],
\end{eqnarray}
which follows easily from 
$$\mu\,e^{2\pi\beta_2 t}\frac \de{\de t}\left (e^{-2\pi\beta_2
t}\widehat{\Gamma_{t,S_r}^\mu}\right )=\hat R_{\beta_2}^\mu
\widehat{\Gamma_{t,S_r}^\mu}.
$$
Choosing $\beta_j=\beta_2$ for $j=2,\dots, k+2,$ we may then take
\eqref{6x} as a substitute for the latter identity in order to perform
the integration by parts argument. Choosing $\beta_2$ appropriately, 
again of the form 
$$\beta_2=\al -\sum_j(2m'_j+1)i\om_j,$$
we find that, except for some trivial remainder terms, 
$K_{\al}^+$ can be written as a finite sum of terms like those in 
\eqref{6w}, only with $\tilde Q_{j_0,k}$ replaced by some polynomial 
in $z''$ of higher degree, and this time with $|q|\ge 2$.

We can proceed in this way, and find that in the $k$-th iteration, the 
$\beta_\ell$'s which have to be chosen are of the form 
\[
\beta_\ell=\al-\sum_j(2m_j+1)i\om_j,
\]
with $\sum_j(2m_j+1)\nu_j<M$ (otherwise, no further integration by parts
 will be needed). This is true for $K_\al^+$. 

In the discussion of $K_\al^-$, one finds in a similar way that the
$\beta_\ell$'s will be of the form $\al+\sum_j(2m_j+1)i\om_j$, again with
$\sum_j(2m_j+1)\nu_j<M$.
Consequently, by Lemma \ref{6C}, $L_{S,\al}$ is locally solvable,
provided 
$L_{S_r}+i(\al\pm\sum_j(2k_j+1)i\om_j)U$ is locally solvable
whenever
$\sum_j(2k_j+1)\nu_j<M$. 

This implies Theorem 2.1 (ii) and completes the
proof of Theorem 2.1.


\setcounter{equation}{0}
\section{Partial results on the case of Hamiltonians with real spectrum}
\label{real}

Our main goal in this section will be to give a proof of Theorem \ref{2E}. Observe that this theorem, in combination with 
Theorem \ref{2A} and Proposition \ref{2B}, immediately also implies Theorem \ref{2G}

Let us thus assume that $S\ne 0 $ has purely real spectrum, that $\Re
Q_S\ge 0, \ N^2=0$ and that property (\ref{C}) is satisfied. Then, by
Corollary \ref{3G}, we have $[S_1,S_2]=0$, hence
$[L_{S_1},L_{S_2}]=-2L_{[S_1,S_2]}U=0$. If we write $\al=a+ib$, with
$a,b$ real, we therefore obtain 
\[
L_S+i\al U=A+iB,
\]
where $A:=L_{S_1}+iaU$ and $B:=L_{S_2}+ibU$ are formally self-adjoint and
commute.

\begin{proposition}\label{7A}
Assume that $A$ and $B$ are left-invariant differential operators on a Lie group $G$ which are formally 
self-adjoint on $C_0^\infty(C_1)\subset L^2(G,dg),$ where $dg$ denotes a
right-invariant Haar measure on $G$. If $A$ and $B$ commute, then $A+iB$
is locally solvable, provided either 
$A$ or $B$ is locally solvable.
\end{proposition}

\noindent{\bf Proof.} It is well-known (compare
\cite{hoermander-grundlehren}, Lemma 6.1.2) that a differential operator
$L$ with smooth coefficients  on $\RR^d$ is locally solvable at $x_0$ if
and only if there exist an open neighborhood $U$ of $x_0$ and constants
$k\in
\NN, C\ge 0$ such that
\begin{equation}\label{7a}
||\vp||_{(-k)} \le C||L^*\vp||_{(k)} \quad \forall \vp\in C_0^\infty (U).
\end{equation}
Here, $||\cdot ||_{(\ell)}$ denotes the Sobolev norm of order $\ell$. If
$L$ is a left-invariant differential operator on $G$, then the
characterization of local solvability given by (\ref{7a}) remains true,
if we assume that $U$ is taken so small that it can be covered by a
single chart, and if we then define the Sobolev norms by means of the
local coordinates, with Lebesgue measure replaced by Haar measure. For
given
$k\in \NN$, we can then find an elliptic right-invariant differential
operator $Q$ on $G$ such that $||\psi||_{(k)}\le ||Q\psi||,\ \psi\in
C_0^\infty(U)$, hence (\ref{7a}) is equivalent to 
\begin{equation}\label{7b}
||\vp||_{(-k)}\le C||Q(L^*\vp)||,\quad  \forall \vp\in
C_0^\infty(U).
\end{equation}
Notice that $QL^*=L^*Q$, since $L^*$ is left-invariant and $Q$ is right-invariant.
Assume now that $A$ and $B$ satisfy the hypotheses of Proposition \ref{7A}, and suppose for instance that $A$ is locally solvable. Then $||A\pm iB)\psi||^2=||A\psi||^2+||B\psi||^2$ for every $\psi\in C_0^\infty(G)$, and we may assume that (\ref{7b}) is satisfied for $L=A$.
This implies
\begin{eqnarray*}
||\vp||_{(-k)}&\le&C||A(Q\vp)||\le C||(A-iB)(Q\vp)||\\
	&=&C||Q((A+iB)^*\vp)||,\ \quad \forall \vp\in C_0^\infty(U),
\end{eqnarray*}
and consequently also $A+iB$ is locally solvable.

\hfill Q.E.D.
\bigskip

\noindent{\bf Proof of Theorem \ref{2E}.}
Consider  $A=L_{S_1}+iaU$. As $S_1^2=0$, by Corollary 3.6, the main
theorem in \cite{mueller-ricci-2} shows that $A$ is locally solvable,
unless
$S_1=0$ and
$a=0$. This implies Theorem \ref{2E}(i), in view of Proposition \ref{7A}.

We have thus reduced ourselves to the case where $A=0,$ i.e. where
$L_S+i\al U=L_{S_2}+ibU.$  But, $L_{S_2}$ is a
real-coefficient operator, and so the remaining cases (ii), (iii) in
Theorem \ref{2E} are immediate consequences of \cite{mueller-ricci-2}. 

\hfill
Q.E.D.

\bigskip

It is perhaps interesting to observe the following corollary to
Proposition \ref{3F}, which in the case  $D=0$ opens up a different
approach to Theorem \ref{2E}.

\begin{proposition}\label{7B}
Assume that $S^2=0$ and $\Re Q_S\ge 0$. Then we can select a symplectic
basis $X_1,\dots,X_n,\ Y_1,\dots,Y_n$ of $V$ such that
$L_S=\sum_{j,k=1}^m b_{jk}Y_jY_k$ for some $m\le n$ and $b_{jk}\in \CC$.
\end{proposition}

\noindent{\bf Proof.} Define the subspaces $W$ and $K$ as in 
Proposition \ref{3F}, and choose a basis $Y_1,\dots,Y_m$ of the isotropic
subspace $W$. Pick $X_1\in V$ such that $\sigma(X_1,Y_j)=\delta_{1,j},
j=1,\dots,m$, and then $X_2$ such that $\sigma(X_2,Y_2)=1$ and $X_2\bot\,
\mbox{span} \{X_1,Y_1,Y_3,\dots,Y_m\}$, and continue in this way to
select $X_j$'s. In the $m$-th step, this means that we pick $X_m$ such
that $\sigma(X_m,Y_m)=1$ and $X_m\bot \,\mbox{span}
\{X_1,\dots,X_{m-1},Y_1,\dots,Y_{m-1}\}$. Then
$U:=\mbox{span}\{X_1,\dots,X_m\}$ is an isotropic subspace too, and $W$
and $U$ are in duality with respect to $\sigma$. In particular, $U\cap
W^\bot=0$, and 
\[
V=U\oplus W\oplus H=U\oplus K,
\]
where $H:=(W\oplus U)^\bot\subset K$ is a symplectic subspace. Moreover,
 since $S_j=N_j, j=1,2$, by the definition of $W$ and $K$ we have
\[
S_j(U)\subset W,\  S_j(K)=0, \quad j=1,2.
\]
Choose a symplectic basis $X_{m+1},\dots,X_n, Y_{m+1},\dots,Y_n$ of $H$. 
Then, with respect to the sets of 
basis elements $X_1,\dots,X_m,Y_1,\dots,Y_m$ and $X_{m+1},\dots,X_n,
Y_{m+1},\dots, Y_n,$ the linear mapping $S=S_1+iS_2$ is represented by a
block matrix of the form
\[
S=\left(
	\begin{array}{cc|cc}
	0 & 0 & 0 & 0\\
	B & 0 & 0 & 0\\
\hline
	0 & 0 & 0 & 0\\
	0 & 0 & 0 & 0
\end{array}
\right),
\]
hence $A=S\left(
	\begin{array}{cc|cc}
	0 & I & 0&0\\
	-I & 0 & 0&0\\
\hline
	0& 0&0 & I\\
	0&0&-I & 0
\end{array}
\right)$ by
$A=
\left(
\begin{array}{cc|cc}
	0 & 0 &0&0\\
	0 & B &0&0\\
\hline
	0 &0& 0&0\\
	0&0&0&0
\end{array}
\right)$.

\noindent This implies the proposition.

\hfill Q.E.D.

\bigskip

Observe that, in suitable coordinates, the operator $L_S+i\al U$, with $S$
 as  in Proposition~\ref{7B}, is a constant coefficient operator, hence
locally solvable, by the Malgrange-Ehrenpreis theorem.
\bigskip

An example which cannot be treated by means of Proposition \ref{7A} is the 
operator $L_S$ from Example \ref{2F}.
Nevertheless, as mentioned in Section \ref{results}, the following
proposition holds true.

\begin{proposition}\label{7C}
Assume that $0< \sqrt{c_1^2+c_2^2}\le m$. Then the operator $L_S+i\al U$,
with $L_S$ given as in Example \ref{2F}, is locally solvable for every
$\al\in \CC$, even though property (\ref{C}) is not satisfied.
\end{proposition}

\noindent{\bf Proof.} The proof will be based on the representation theoretic 
criterion given by Theorem 4.1 in \cite{mueller-peloso-nonsolv}.
Adapting the notation from \cite{mueller-peloso-nonsolv}, denote for 
$\mu\in
\RR\setminus \{0\}$ by
$\pi_\mu$ the Schr\"odinger representation of $\HH_2$ on $L^2(\RR^2)$,
with parameter $\mu$. For its differential, we have
\[
d\pi_\mu(X_j)=\frac \de{\de x_j},\ d\pi_\mu(Y_j)=i\mu x_j,\ (j=1,2), \quad
d\pi_\mu(U)=i\mu,
\]
if we denote the coordinates in $\RR^2$ by $x=(x_1,x_2)$. Putting 
$\L:=L_S+i\al U$, then 
\begin{eqnarray}\label{7c}
 &&-d\pi_\mu(\L)\\
&&\quad=\mu^2((m+c_1)x_1^2+(m-c_1)x_2^2+2c_2x_1x_2)
+2\mu(x_2\de _{x_1}-x_1\de_{x_2})-\al\mu.\nonumber
\end{eqnarray}
Introducing polar coordinates $x_1=r\cos \theta, x_2=r\sin \theta,\  r\ge
0$, we obtain
\begin{equation}\label{7d}
d\pi_\mu(\L)=\mu[2\de _\theta +\al-\mu r^2q(\theta)],
\end{equation}
where $q(\theta):=m+c_1\cos (2\theta)+c_2\sin (2\theta)$.
Observe that $q\ge 0$, since $m\ge \sqrt{c_1^2+c_2^2}$. 

\noindent For $\theta_0$ fixed, we put
$$Q_{\theta_0}(\theta):=\int\limits_{\theta_0}^\theta
q(\theta)\, d\theta=m(\theta-\theta_0)+\psi(\theta-\psi(\theta_0),$$ where
$\psi(\theta):=\frac {c_1}2 \sin (2\theta)-\frac{c_2}2\cos (2\theta)$.
Notice that 
\begin{equation}\label{7e}
Q_{\theta_0}(\theta)\ge 0, \quad\mbox{ if } \theta\ge \theta_0,
\end{equation}
and that $\psi$ is $2\pi$-periodic. Define
$F_{\theta_0}^\mu(r,\theta)$ to be the $2\pi$-periodic extension of
$\tilde F_{\theta_0}^\mu(r,\theta)$ given by 
\begin{equation}\label{7f}
\tilde F_{\theta_0}^\mu(r,\theta):=
e^{-\frac \al 2(\theta-\theta_0)}e^{\frac 1 2\mu
r^2Q_{\theta_0}(\theta)},\quad
\theta_0\le \theta < \theta_0+2\pi.
\end{equation}
Then one checks easily that, in the sense of ($2\pi$-periodic) distributions,
\begin{eqnarray*}
\lefteqn{d\pi\mu(\L)[h(r)F_{\theta_0}^\mu(r,\theta)]}\\
&&= 2\mu h(r)[F_{\theta_0}^\mu(r,\theta_0)-F_{\theta_0}^\mu(r,\theta_0+2\pi)]\delta_{\theta_0}(\theta)\\
&&=2\mu h(r)\ [1-e^{-\pi\al+\pi\mu r^2m}]\delta_{\theta_0}(\theta).
\end{eqnarray*}
Notice that for $v=(r_0\cos \theta_0,r_0\sin \theta_0)\ne 0$, we have
$\delta_v(x)=\frac 1{r_0}\delta_{r_0}(r)\delta_{\theta_0}(\theta)$. 
Therefore, if we put
\begin{equation}\label{7g}
H_\mu(x,v):=\frac{F_{\theta_0}^\mu (r_0,\theta)}
	{2\mu r_0[1-e^{-\pi\al+\pi\mu r_0^2m}]} \delta_{r_0}(r),
\end{equation}
we have $H_\mu(\cdot,v)\in \S'(\RR^2)$, for a.e. $v$, and
$d\pi_\mu(\L)H_\mu(\cdot,v)=\delta(v-\cdot)$. Write
$H_\mu=G_\mu\delta_{r_0}(r)$,  and $\al=a+ib,\ a,b\in \RR$. We claim that 
\begin{equation}\label{7h}
G_\mu(x,v)=\left\{
	\begin{array}{ll}
	\frac{g_\mu(x,v)}{\mu r_0},& \mbox{ if $b\not\in 2\ZZ$,}\\
	\frac{g_\mu(x,v)}{\mu r_0(a-\mu mr_0^2)},& \mbox{ if $b\in 2\ZZ$,}
\end{array}
\right.
\end{equation}
where 
\begin{equation}\label{7i}
||g_\mu||_\infty \le C_\al,
\end{equation}
with $C_\al$ independent of $\mu$.

Indeed, if $\mu<0$, then $|F_{\theta_0}^\mu|\le C_\al$, since
$Q_{\theta_0}\ge 0$ by (\ref{7e}). Moreover, if $b\not\in 2\ZZ$, then 
$|1-e^{-\pi\al+\pi\mu r_0^2m}|$ is bounded from below, and if $b\in 2\ZZ$,
then $1-e^{-\pi\al+\pi\mu r_0^2m}=1-e^{-\pi(a-\mu mr_0^2)}$, and
(\ref{7h}), (\ref{7i}) follow immediately.

On the other  hand, if $\mu>0$, then for $\mu r_0^2m>>1$ we have
\[
|\mu r_0G_\mu(x,v)|\le C_\al e^{\frac \mu 2 r_0^2Q_{\theta_0}(\theta)
-\pi\mu r_0^2m}.
\]
Moreover, $Q_{\theta_0}(\theta)\le Q_{\theta_0}(\theta_0+2\pi)=2\pi m$,
and thus $|\mu r_0 G_\mu(x,v)|\le C_\al$. Finally, the case where 
$\mu r_0^2 m\le C$ is easy, and again we obtain (\ref{7h}), (\ref{7i}).

(\ref{7h}) shows that singularities may arize, if $b\in 2\ZZ$, and even
 if $b\not\in 2\ZZ$, it turns out that some negative powers of $\mu$ may
arize in a later step of the proof. For given $N\in \NN$, we therefore
define
\[
G_\mu ^N(x,v):=\left\{
\begin{array}{ll}
(i\mu)^NG_\mu(x,v), & \mbox{ if $b\not\in 2\ZZ$,}\\
(i\mu)^N(i\mu a-i\mu^2mr_0^2), & \mbox{ if $b\in 2\ZZ$,}
\end{array}
\right.
\]
and $H_\mu ^N=G_\mu^N \delta_{r_0}(r)$. Then
\begin{eqnarray}\label{7j}
d\pi_\mu(\L)H_\mu^N(\cdot,v)&=&(i\mu)^N(i\mu a-i\mu^2 m|v|^2)^\sigma \delta(v-\cdot)\\
&=& d\pi _\mu (U^N(aU+im(Y_1^2+Y_2^2))^\sigma)\delta(v-\cdot),\nonumber
\end{eqnarray}
where $\sigma:=0$, if $b\not\in 2\ZZ$, and $\sigma=1$, if $b\in 2\ZZ$.
For $\vp\in C_0^\infty (\RR^2\times \RR^2)$, it then follows easily from (\ref{7h}) that the mapping $v\mapsto 
\langle H_\mu^N(\cdot,v),\vp(\cdot,v)\rangle$, which, in polar coordinates, is given by
\[
\int_0^{2\pi} G_\mu^N((r_0,\theta),(r_0,\theta_0))
\,\vp((r_0,\theta),(r_0,\theta_0))\, r_0\, d\theta,
\]
is integrable with respect to $v$ and defines a tempered distribution
$H_\mu^N,$ namely
 \begin{eqnarray*} &&\l H_\mu^N,\vp\r\\
 &&:=\int \l H_\mu^N(\cdot,v),\vp(\cdot,v)\r dv
=\int_0^{+\infty}\int_0^{2\pi}\int_0^{2\pi} 
G_\mu^N((r_0,\theta),(r_0,\theta_0))\vp((r_0,\theta))\,d\theta \,d\theta_0
r_0^2\,dr,
\end{eqnarray*}
 for every $\mu\ne 0$. We set $\tilde
H_\mu^N(x,\eta):=H_\mu^N(\frac \eta\mu+\frac x 2, \frac \eta \mu-\frac x
2)$, where the change of variables is to be understood in the sense of
distributions, i.e.
\[
\l\tilde H_\mu^N,\vp\r=\l H_\mu^N,\vp(u-v,\frac\mu 2(u+v))\r|\mu|^2.
\]
Observe that, by (\ref{7h}), (\ref{7i}), we can find a continuous Schwartz 
norm $||\cdot ||_\S$ on $\S(\RR^2\times \RR^2)$ and $M\in \NN$ such that
\[
|\l H_\mu^N,\psi\r|\le |\mu|^N(|\mu|^M+|\mu|^{-M})||\psi||_\S,\quad \psi
\in
\S,
\]
for every $N$. This shows that
\[
\l \tilde H^N,f\r:=\int\limits_{\RR\setminus \{0\}}\l \tilde
H_\mu^N,\psi(\cdot,\mu)\r\, d\mu,\quad f\in \S(\RR^2\times \RR^2\times
\RR)
\]
defines a tempered distribution on $\RR^2\times \RR^2\times \RR$, provided we choose $N$ sufficiently large.
But then the $H_\mu^N$ satisfy the hypothesis of
\cite{mueller-peloso-nonsolv}, Theorem 4.1, except for condition (i) in
this theorem, which has to be replaced by (\ref{7j}). The proof of
Theorem 4.1 in \cite{mueller-peloso-nonsolv} then still applies and shows
that there is a tempered distribution $K^N\in
\S'(\HH_2)$ such that 
\[
\L K^N=cU^N(aU+im(Y_1^2+Y_2^2))^\sigma \delta_0,
\]
for a suitable constant $c\ne 0$.

Since the operator $U^N(aU+im(Y_1^2+Y_2^2))^\sigma$ is locally solvable, 
this implies local solvability of $\L$ (compare e.g.
 the proof of Lemma 7.4 in \cite{mueller-ricci-cone}).

\hfill Q.E.D.

\bibliographystyle{plain}

\vskip1cm

\noindent\emph{Mathematisches Seminar,  C.A.-Universit\"at Kiel,
Ludewig-Meyn-Str.4, D-24098 Kiel, Germany\\
 e-mail: mueller@math.uni-kiel.de}
\end{document}